\RequirePackage{fix-cm}
\documentclass[smallextended]{svjour3}       % onecolumn (second format)
\smartqed  % flush right qed marks, e.g. at end of proof
\usepackage{graphicx}
\usepackage[utf8]{inputenc}
\usepackage{amsmath}
\usepackage[english]{babel}
\usepackage{amssymb}
\usepackage{amsfonts}
\usepackage{bm}
\usepackage{graphicx}
\usepackage{xcolor}
\usepackage[hidelinks=true]{hyperref}
\usepackage{parskip}
\usepackage{multicol}
\usepackage{float}
\usepackage{soul}

\usepackage{verbatim}

\floatstyle{plaintop}
\restylefloat{table}
\usepackage[tableposition=top]{caption}

\usepackage{chngcntr}
\usepackage{apptools}

\AtAppendix{\counterwithin{proposition}{subsection}}
\AtAppendix{\counterwithin{equation}{subsection}}

\usepackage{booktabs}

\usepackage{natbib}
\hypersetup{urlcolor=blue, citecolor=blue, colorlinks=true, linkcolor=blue} 

\numberwithin{equation}{section}
\numberwithin{figure}{section}
\numberwithin{table}{section}

\usepackage{blindtext}
\usepackage{lineno}

\begin{document}

%\title{The Turning Arcs Algorithm to Simulate Isotropic Gaussian Random Fields on Spheres and Hyperspheres}

\title{The Turning Arcs: a Computationally Efficient Algorithm to Simulate Isotropic Vector-Valued Gaussian Random Fields on the $d$-Sphere}

\authorrunning{The Turning Arcs Simulation Algorithm}  

\author{Alfredo Alegr\'ia        \and
        Xavier Emery \and  Christian Lantu\'ejoul %etc.
}

\institute{A. Alegr\'ia (Corresponding author)  \at
             Departamento de Matem\'atica, 
              Universidad T{\'e}cnica Federico Santa Mar{\'i}a, Valpara{\'i}so, Chile
              \email{alfredo.alegria@usm.cl}           %  \\
           \and
           X. Emery  \at 
           Department of Mining Engineering, University of Chile, Santiago, Chile\\
            Advanced Mining Technology Center, University of Chile, Santiago, Chile
                            \and
           C. Lantu\'ejoul  \at
              Centre de G\'eosciences, MINES ParisTech, PSL University, Paris, France
}

\date{Received: date / Accepted: date}

\maketitle

\begin{abstract}
Random fields on the sphere play a fundamental role in the natural sciences. This paper presents a simulation algorithm parenthetical to the spectral turning bands method used in Euclidean spaces, for simulating scalar- or vector-valued Gaussian random fields on the $d$-dimensional unit sphere. The simulated random field is obtained by a sum of Gegenbauer waves, each of which is variable along a randomly oriented arc and constant along the parallels orthogonal to the arc. Convergence criteria based on the Berry-Ess\'{e}en inequality are proposed to choose suitable parameters for the implementation of the algorithm, which is illustrated through numerical experiments. A by-product of this work is a closed-form expression of the Schoenberg coefficients associated with the Chentsov and exponential covariance models on spheres of dimensions greater than or equal to $2$.

\keywords{Schoenberg sequence \and Turning Bands \and Gegenbauer polynomials \and Central limit approximation \and Berry-Ess\'{e}en inequality}

\end{abstract}

\section{Introduction}
 
Spherically indexed Gaussian random fields have attracted a growing interest in recent decades. They are useful in the modeling of georeferenced variables arising in many branches of applied sciences, such as astronomy, climatology, oceanography, biology and geosciences, amongst many others. We refer the reader to \cite{marinucci2011random}, \cite{jeong2017spherical} and \cite{porcu2018modeling} for recent reviews about this topic. In general, the space consists of a $2$-dimensional sphere, but hyperspheres are sometimes met, e.g., in high-dimensional shape analysis \citep{Dryden, Mardia}.

Simulation is crucial for the development of new applications in spatial statistics. It is well known that simulation algorithms based on the Cholesky decomposition of the covariance matrix \citep{Ripley} are computationally prohibitive when the sample size is large, since the order of computation of the Cholesky decomposition is equal to the cube of the sample size. As a result, the search for new efficient methods to simulate Gaussian random fields in spherical domains is of paramount importance.  Within the class of isotropic random fields, i.e., random fields whose finite-dimensional distributions are invariant under rotations, several appealing alternatives have been proposed, including spherical harmonic approximations \citep{lang2015isotropic, de2018regularity, emery2019, lantuejoul2019}, circulant embedding approaches \citep{cuevas2019fast}, random coin type methods \citep{hansen2015gaussian}, and simulations over Euclidean spaces restricted to low-dimensional spheres \citep{emery2019turning}.  

In this paper, we propose a simple algorithm that simulates a Gaussian random field with a prescribed isotropic covariance structure, based on adequate combination of Gegenbauer waves. Our proposal, named the `turning arcs' method, can be seen as the spherical counterpart of the spectral turning bands method developed in Euclidean spaces (see, e.g.,  \citealp{matheron1973intrinsic}, \citealp{mantoglou1982turning}, \citealp{lantu2002}, \citealp{emery2006} and \citealp{emery2016}). The advantages of this algorithm over existing ones are threefold:
\begin{itemize}
\item[1.] It is computationally less expensive than approximations based on spherical harmonics.
\item[2.] It is applicable to the simulation not only on the $2$-sphere, but also on the $d$-sphere, for any dimension $d$.
\item[3.] It allows the simulation not only of scalar random fields, but also on vector random fields.
\end{itemize}
% \fbox{advantages of our proposal with respect to previous alternatives}

The outline of the paper is as follows. In Section \ref{rfs} preliminary results about isotropic scalar- and vector-valued Gaussian random fields on the $d$-sphere are reviewed. The `turning arcs' simulation algorithm is then presented in Section \ref{section-sim}. In Section \ref{examples}, the applicability of our proposal is illustrated through numerical examples. Section \ref{practicalaspects} discusses the computational implementation and provides some  guidelines to practitioners. Section \ref{conclusions} concludes the paper, while technical proofs are given in Appendices.

\section{Background}

\label{rfs}

\subsection{Scalar-Valued Isotropic Gaussian Random Fields on the Sphere}

Let $\mathbb{S}^d  = \{\bm{x} \in \mathbb{R}^{d+1}: \bm{x}^\top \bm{x}=1\}$ be the $d$-dimensional unit sphere embedded in $\mathbb{R}^{d+1}$, where $^\top$ denotes the transpose operator, and consider a real-valued random field, $Z = \{ Z(\bm{x}): \bm{x}\in\mathbb{S}^d\}$ with finite second-order moments.  We assume  $Z$ to be Gaussian, i.e., for all $k\in\mathbb{N}^*$ and $\bm{x}_1,\hdots,\bm{x}_k\in\mathbb{S}^d$, the random vector $\{Z(\bm{x}_1),\hdots,Z(\bm{x}_k)\}^\top$ follows a multivariate Gaussian distribution. Thus, $Z$  is completely characterized by its mean function and its covariance function given by $$C(\bm{x}_1,\bm{x}_2) = {\rm cov}\{Z(\bm{x}_1),Z(\bm{x}_2)\}, \qquad \bm{x}_1,\bm{x}_2\in\mathbb{S}^d.$$ 
Let us introduce the geodesic distance on $\mathbb{S}^d$, which is the main ingredient to define the property of isotropy of a random field. For two locations, $\bm{x}_1$ and $\bm{x}_2$ in $\mathbb{S}^d$,  their geodesic distance is defined as $\vartheta(\bm{x}_1,\bm{x}_2) = \arccos\{ \bm{x}_1^\top \bm{x}_2 \} \in [0,\pi]$.  We shall equivalently use $\vartheta(\bm{x}_1,\bm{x}_2)$ or the shortcut $\vartheta$ to denote the geodesic distance. Following \citet{marinucci2011random}, the random field is called (weakly) isotropic if it has constant mean and if its covariance function can be written as
\begin{equation}
\label{isotropic-part}
C(\bm{x}_1, \bm{x}_2) = K\{ \vartheta(\bm{x}_1,\bm{x}_2) \}, \qquad \bm{x}_1, \bm{x}_2 \in \mathbb{S}^d,
\end{equation} 
for some continuous function $K: [0,\pi]\rightarrow \mathbb{R}$. Thus, the covariance function just depends on the geodesic distance. It is common to call $K$ the \emph{isotropic part} of the covariance function $C$ (see, e.g., \citealp{guella2018unitarily}). For Gaussian random fields, isotropy also implies that the probability distribution of $\{Z(\bm{x}_1),\hdots,Z(\bm{x}_k)\}^\top$ is invariant under the group of rotations on $\mathbb{S}^d$ (see \citealp{marinucci2011random}).    

Positive semi-definiteness is a necessary and sufficient condition for a function to be a valid covariance.
%Any candidate function to be a covariance function must be positive semi-definite, which is a necessary and sufficient condition to be a valid covariance.
% technical condition. Indeed, note that for all $k\in\mathbb{N}$, and for all systems of points $\bm{x}_1,\hdots, \bm{x}_k \in \mathbb{S}^d$ and constants $a_1,\hdots,a_k \in \mathbb{R}$, we have that
%\begin{equation}
%\label{pos_def}
% {\rm var}\left\{  \sum_{i=1}^k a_i Z(\bm{x}_i)\right\} = \sum_{i=1}^k \sum_{j=1}^k a_i a_j C(\bm{x}_i,\bm{x}_j) \geq 0.
%\end{equation}
% Condition (\ref{pos_def}), called \emph{semi positive definiteness},  is a necessary and sufficient condition for a covariance function.  
In his pioneering paper,  \citet{schoenberg1942} showed that $C$ as in (\ref{isotropic-part}) is  positive semi-definite if, and only if, its isotropic part $K$ has a series representation of the form
\begin{equation}
\label{schoenberg1}
%K(\vartheta) = \sum_{n=0}^\infty   b_{n,d}   \frac{{G}_n^{(d-1)/2}(\cos \vartheta)}{ {G}_n^{(d-1)/2}(1)}, \qquad   0 \leq \vartheta  \leq \pi,
K(\vartheta) = \sum_{n=0}^\infty   b_{n,d} \,  {G}_n^{(d-1)/2}(\cos \vartheta), \qquad   0 \leq \vartheta  \leq \pi,
\end{equation}
where $\{b_{n,d}: n\in\mathbb{N}\}$ is a sequence of nonnegative coefficients such that $\sum_{n=0}^\infty b_{n,d} {G}_n^{(d-1)/2}(1) < +\infty$, referred to as a \emph{Schoenberg sequence} \citep{gneiting2013strictly}, while $\{{G}_n^{\lambda}: n\in\mathbb{N}\}$ is the sequence of $\lambda$-Gegenbauer polynomials \citep{abramowitz}, which are implicitly defined through the identity
$$ \frac{1}{(1-2rt + t^2)^\lambda}  = \sum_{n=0}^\infty G_n^{\lambda}(r)t^n, \qquad -1\leq r \leq 1. $$
%Here, $\{b_{n,d}: n\in\mathbb{N}_0\}$ is a sequence of nonnegative coefficients, such that $\sum_{n=0}^\infty b_{n,d} < \infty$. As in \cite{gneiting2013strictly}, such a sequence is referred to as a \emph{Schoenberg sequence}.
%we refer to the sequence $\{b_{n,d}: n\in\mathbb{N}_0\}$ as a \emph{Schoenberg sequence}. 

The Gegenbauer polynomials can be calculated in a straightforward manner by use of the following recurrence relationships:
\begin{equation}
    \label{Gegenbauer}
    \begin{cases}
    & {G}_0^{\lambda}(r) = 1; \\
    & {G}_1^{\lambda}(r) = 2 \, \lambda \, r;\\
    & {G}_n^{\lambda}(r) = \frac{2(n+\lambda-1)}{n} \, r \, \, {G}_{n-1}^{\lambda}(r) - \frac{n+2\lambda-2}{n} \, {G}_{n-2}^{\lambda}(r), \qquad n>1.
    \end{cases}
\end{equation}

%Also, the following inequality holds %\citep{abramowitz}:
%\begin{equation}
%    \label{Gegenbauer2}
%\lvert {G}_n^{\lambda}(r) \rvert \leq {G}_n^{\lambda}(1) = \frac{\Gamma(2\lambda+n)}{\Gamma(2\lambda) \Gamma(n+1)}, \qquad -1\leq r \leq 1, \, n\in\mathbb{N}.
%\end{equation}

In practice, the most usual cases correspond to spheres of dimensions $d=1$ or $d=2$. When $d=1$, Schoenberg's expansion is written in terms of Chebyshev polynomials, $G_n^0(\cos \vartheta) = \cos(n \vartheta)$. When $d=2$, one obtains an expansion in terms of Legendre polynomials, $G_n^{1/2}(\cos \vartheta) = P_n(\cos \vartheta)$.  \\
 
 There is a one-to-one correspondence between an isotropic covariance $K$ and its Schoenberg sequence. Classical inversion formulae yield the identity \citep{schoenberg1942, gneiting2013strictly, ziegel2014convolution}
 %\begin{equation}
 %b_{n,d} = \frac{(2n+d-1)}{2^{3-d}\pi} \frac{\Gamma((d-1)/2)^2}{\Gamma(d-1) G_n^{(d-1)/2}(1)} \int_0^\pi G_n^{(d-1)/2}(\cos \vartheta) (\sin \vartheta )^{d-1} K(\vartheta) \text{d}\vartheta, \qquad n\in\mathbb{N},
 %\end{equation}
 \begin{equation}
 \label{bnd1}
 b_{n,d} = \frac{1}{\parallel G_n^{(d-1)/2} \parallel^2} \int_0^\pi G_n^{(d-1)/2}(\cos \vartheta) (\sin \vartheta )^{d-1} K(\vartheta) \text{d}\vartheta, \qquad n\in\mathbb{N},
 \end{equation}\\
 with (\citet{abramowitz}, formula 22.2.3)
\begin{equation}
\label{bnd2}
\begin{split}
\parallel G_n^{(d-1)/2} \parallel^2 &= \int_0^\pi \bigl[ G_n^{(d-1)/2} (\cos \vartheta) \bigr]^2 \, \sin^{d-1} (\vartheta) \, d \vartheta \\
&= \begin{cases}
\frac{2\pi}{n^2} &\text{if $d=1$} \\
\frac{2^{3-d}\pi}{(2n+d-1)} \frac{\Gamma(d-1+n)}{n! \Gamma((d-1)/2)^2} &\text{if $d \geq 2$.}
\end{cases}
\end{split}
\end{equation}

%For $d=1$, we have that $$ b_{0,1} = \frac{1}{\pi} \int_0^\pi  K(\vartheta)  \text{d}\vartheta  \qquad \text{ and } \qquad  b_{n,1} = \frac{2}{\pi} \int_0^\pi  K(\vartheta) \cos(n\theta)  \text{d}\vartheta,$$
%for all  $n\in\mathbb{N}$. 

The isotropic covariance function, or its Schoenberg sequence, is often specified to belong to a parametric family whose members are known to be positive semi-definite.  For a thorough review on positive semi-definite functions on spheres and a list of parametric families, we refer the reader to \citet{hagzhgrbn11}, \citet{gneiting2013strictly}, \citet{artggipcu18} and \citet{lantuejoul2019}. The Schoenberg sequences of two specific parametric families (Chentsov and exponential covariances) on $\mathbb{S}^d$, $d \geq 2$, are also given in Appendix \ref{sub:sbg}, which seems to be a new result.

\subsection{Vector-Valued Isotropic Gaussian Random Fields on the Sphere}

We now turn to a description of vector-valued random fields. Let $\bm{Z}=\{ [Z_1(\bm{x}), \hdots, Z_p(\bm{x})]^\top: \bm{x}\in\mathbb{S}^d\}$ be a $p$-variate random field, with each component having finite second-order moments. 
We assume $\bm{Z}$ to be Gaussian, i.e., for all $k\in\mathbb{N}^{*}$ and $\bm{x}_1,\hdots,\bm{x}_k\in\mathbb{S}^d$, the random vector $\{\bm{Z}(\bm{x}_1),\hdots,\bm{Z}(\bm{x}_k)\}^\top$ follows a multivariate Gaussian distribution, where $\bm{Z}(\bm{x}) =  [Z_1(\bm{x}), \hdots, Z_p(\bm{x})]^\top$. 
%We assume $\bm{Z}(\bm{x}) =  [Z_1(\bm{x}), \hdots, Z_p(\bm{x})]^\top$ to be a  Gaussian random field. Let $C_{ij}(\bm{x}_1,\bm{x}_2) = {\rm cov}\{Z_i(\bm{x}_1),Z_j(\bm{x}_2)\}$,  for all $\bm{x}_1,\bm{x}_2\in\mathbb{S}^d$ and $i,j=1,\hdots,p$. 
We denote by $\bm{C}(\bm{x}_1,\bm{x}_2)$ the $p\times p$ covariance matrix between $\bm{Z}(\bm{x}_1)$ and $\bm{Z}(\bm{x}_2)$, with the $(i,j)$th entry equal to $C_{ij}(\bm{x}_1,\bm{x}_2)$. The diagonal elements, $C_{ii}(\bm{x}_1,\bm{x}_2)$, are called direct covariance functions, whereas the off-diagonal elements, $C_{ij}(\bm{x}_1,\bm{x}_2)$, for $i\neq j$, are called cross-covariance functions.

The isotropy of a vector-valued random field can be defined in a similar fashion to the scalar-valued case. Indeed, the random field $\bm{Z}$  is called isotropic if each of its components has a constant mean and if its matrix-valued covariance function can be written as
\begin{equation*}            
\label{isotropic-part-vector}
\bm{C}(\bm{x}_1, \bm{x}_2) = \bm{K}\{ \vartheta(\bm{x}_1,\bm{x}_2) \}, \qquad \bm{x}_1, \bm{x}_2 \in \mathbb{S}^d,
\end{equation*} 
for some continuous matrix-valued function $\bm{K}: [0,\pi]\rightarrow \mathbb{R}^{p\times p}$. The condition of positive semi-definiteness can also be adapted to the vector-valued case.
% as follows. For all $k\in\mathbb{N}$, and for all systems of points $\bm{x}_1,\hdots, \bm{x}_k \in \mathbb{S}^d$ and vectors $\bm{a}_1,\hdots, \bm{a}_k \in \mathbb{R}^p$, we have that
%\begin{equation*}
%\label{pos_def_vector}
% {\rm var}\left\{  \sum_{i=1}^k \bm{a}_i^\top \bm{Z}(\bm{x}_i)\right\} = \sum_{i=1}^k \sum_{j=1}^k \bm{a}_i^\top  \bm{C}(\bm{x}_i ,\bm{x}_j)  \bm{a}_j \geq 0.
%\end{equation*}
The Schoenberg's expansion for the matrix-valued isotropic part is given by \citep{Yaglom:1987,hannan2009multiple}
\begin{equation}
\label{schoenberg1_vector}
%\bm{K}(\vartheta) = \sum_{n=0}^\infty   \bm{B}_{n,d}   \frac{{G}_n^{(d-1)/2}(\cos \vartheta)} {{G}_n^{(d-1)/2}1)}, \qquad   0 \leq \vartheta \leq \pi,
\bm{K}(\vartheta) = \sum_{n=0}^\infty   \bm{B}_{n,d} \,  {G}_n^{(d-1)/2}(\cos \vartheta), \qquad   0 \leq \vartheta \leq \pi,
\end{equation}
where $\{\bm{B}_{n,d}: n\in\mathbb{N}\}$ is a sequence of positive semi-definite matrices (called Schoenberg matrices) such that $\sum_{n=0}^\infty \bm{B}_{n,d} {G}_n^{(d-1)/2}(1) < +\infty$ (element-wise summation). Similarly to the scalar-valued scenario, Fourier calculus implies that 
 $$  \bm{B}_{n,d} = \frac{1}{\parallel G_n^{(d-1)/2} \parallel^2}\int_0^\pi G_n^{(d-1)/2}(\cos \vartheta) (\sin \vartheta )^{d-1} \bm{K}(\vartheta) \text{d}\vartheta, \qquad n\in\mathbb{N}.$$
%When $d=1$, we have the inversion formula $$ \bm{B}_{0,1} = \frac{1}{\pi} \int_0^\pi  \bm{K}(\vartheta)  \text{d}\vartheta  \qquad \text{ and } \qquad  \bm{B}_{n,1} = \frac{2}{\pi} \int_0^\pi  \bm{K}(\vartheta) \cos(n\theta)  \text{d}\vartheta,$$
%for all  $n\in\mathbb{N}$. 

\section{The Turning Arcs Simulation Algorithm}
\label{section-sim}

\subsection{Scalar-Valued Case}

This section presents an algorithm for simulating scalar-valued isotropic Gaussian random fields on $\mathbb{S}^d$. The representation (\ref{schoenberg1}) allows for an immediate simulation procedure based on the Schoenberg sequence $\{b_{n,d}: n\in \mathbb{N}\}$. The following proposition is crucial to develop the simulation algorithm. 

\begin{proposition}
\label{prop1}
 Let $\varepsilon$ be a random variable with zero mean and unit variance, $\bm{\omega}$ a random vector uniformly distributed on $\mathbb{S}^d$, and $\kappa$  a discrete random variable with $\mathbb{P}(\kappa = n) = a_n$,  $n\in\mathbb{N}$, where $\mathbb{P}$ indicates the probability. Suppose that the support of the probability mass sequence $\{a_n: n\in\mathbb{N}\}$ contains the support of the Schoenberg sequence $\{b_{n,d}: n\in\mathbb{N}\}$ and that $\varepsilon$, $\bm{\omega}$ and $\kappa$ are independent. Then, 
\begin{enumerate}
\item[(1)]  For $d=1$, the random field defined by
\begin{equation}
\label{simulacion_d=1}
{Z}(\bm{x}) = \varepsilon \sqrt{ \frac{2 b_{\kappa,1}}{a_\kappa } }    \cos(\kappa \vartheta(\bm{\omega}, \bm{x})), \qquad \bm{x}\in\mathbb{S}^1,
\end{equation}
is isotropic, with zero mean 
%$\mathbb{E}\{ {Z}(\bm{x})  \} = \sqrt{2 a_0 b_{0,1}}$, for all $\bm{x}\in\mathbb{S}^1$, 
and covariance function with isotropic part given by
$$K(\vartheta) =  \sum_{n=0}^\infty b_{n,1} \cos(n \vartheta ),    \qquad 0\leq \vartheta \leq \pi.$$
%$$K(\vartheta) =  2 (1-a_0) b_{0,1}  + \sum_{n=1}^\infty b_{n,1} \cos(n \vartheta ),    \qquad 0\leq \vartheta \leq \pi.$$

\item[(2)] For $d\geq 2$, the random field defined by
\begin{equation}
\label{simulacion2}
{Z}(\bm{x}) = \varepsilon \sqrt{\frac{b_{\kappa,d} (2\kappa + d -1)}{a_\kappa (d-1) } }   {G}_\kappa^{(d-1)/2}(\bm{\omega}^\top \bm{x}), \qquad \bm{x}\in\mathbb{S}^d,
\end{equation}
is isotropic, with zero mean % $\mathbb{E}\{ {Z}(\bm{x})  \} = \sqrt{a_0 b_{0,d}}$, for all $\bm{x}\in\mathbb{S}^d$, 
and covariance function with isotropic part given by (\ref{schoenberg1}).
%$$ K(\vartheta) =  \sum_{n=0}^\infty b_{n,d} {G}_n^{(d-1)/2}( \cos \vartheta),     \qquad 0\leq \vartheta \leq \pi.$$
%$$ K(\vartheta) =  (1-a_0) b_{0,d}  + \sum_{n=1}^\infty b_{n,d} {G}_n^{(d-1)/2}( \cos \vartheta),     \qquad 0\leq \vartheta \leq \pi.$$
\end{enumerate} 
\end{proposition}

Proposition \ref{prop1}, the proof of which is deferred to Appendix \ref{proof-prop1} for a neater exposition, provides a procedure to simulate isotropic random fields on the sphere with the predefined covariance function (\ref{schoenberg1}). 
%, except for a slight modification in the first coefficient of its expansion that does not affect the flexibility of this simulation strategy. 
Note that the algorithm separates the choice of the adaptive Schoenberg sequence, which provides the covariance structure of the simulated random field, from the choice of the probability mass sequence $\{a_n: n\in\mathbb{N}\}$ according to which the degrees of the Gegenbauer polynomials are simulated. %, which is equivalent to an importance sampling technique. %We observe that location and dispersion parameters can also be added in order to control the mean and variance of the simulated random field. 

The simulated random field reproduces the desired first- and second-order moments (zero mean and isotropic covariance $K$), but is not normally distributed. A central limit approximation of a Gaussian random field with the same first- and second-order moments can be obtained by \citep{lantu2002,Chiles}
\begin{equation}
\label{clt1}
\widetilde{Z}(\bm{x}) =   \frac{1}{\sqrt{L}}\sum_{\ell = 1}^L {Z}_\ell(\bm{x}), \qquad \bm{x}\in\mathbb{S}^d,
\end{equation}
where $L$ is a large integer and ${Z}_1(\bm{x}),\hdots, {Z}_L(\bm{x})$ are $L$ independent copies of ${Z}(\bm{x})$. 

\begin{figure}[htp]
    \begin{center}
    \includegraphics[width=10cm]{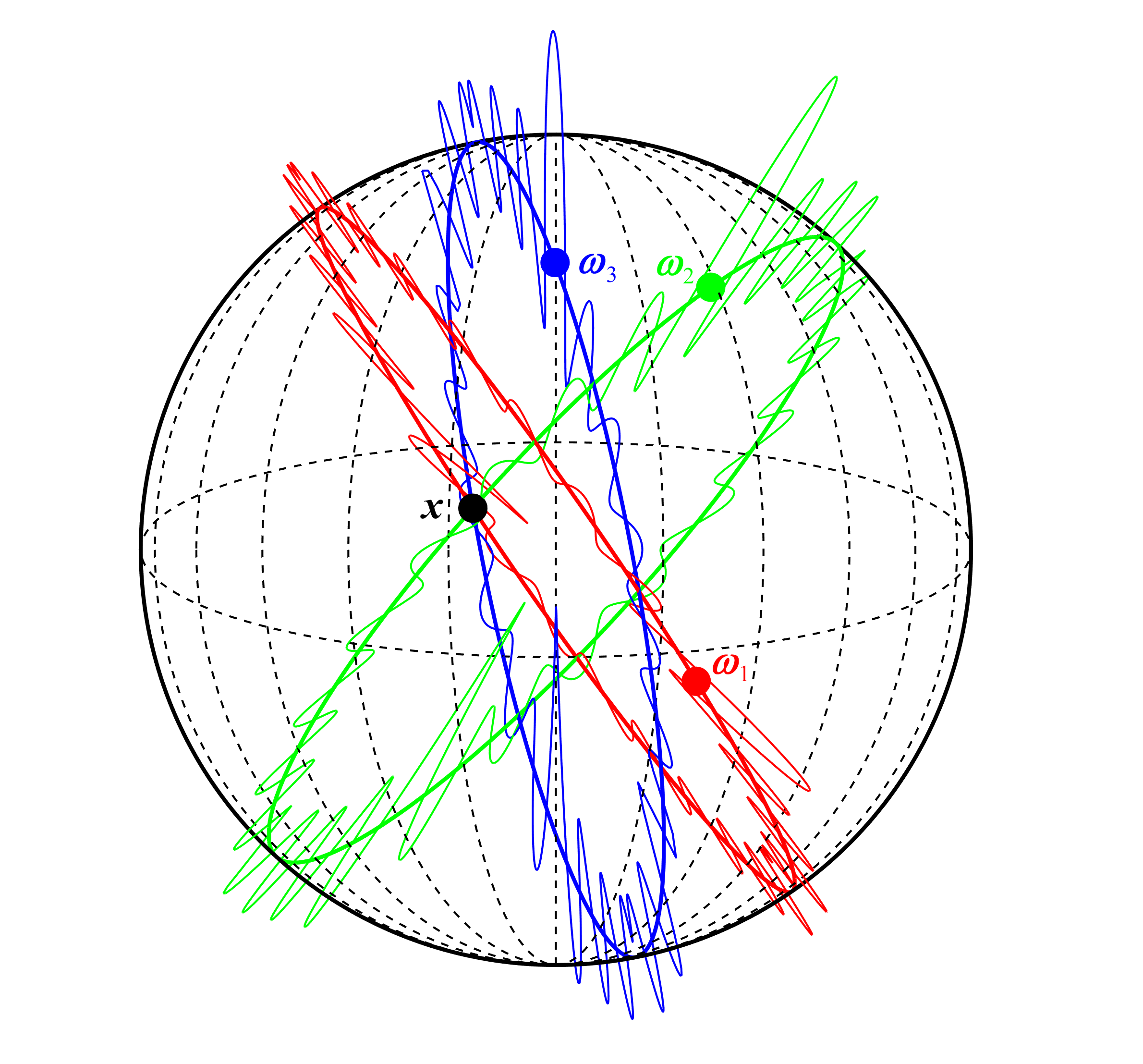}
    \end{center}
\caption{Turning arcs on the $2$-sphere: three arcs with random poles $\bm{\omega}_1$, $\bm{\omega}_2$ and $\bm{\omega}_3$ passing through a point $\bm{x}$ (red, green and blue great circles) and the basic random fields ${Z}_1$, ${Z}_2$ and ${Z}_3$ (thin colored lines) varying along these arcs. The equator and a few meridians are superimposed (dashed lines). The simulated random field ${Z}$ at $\bm{x}$ is a weighted sum of the three basic random fields at this point} 
\label{fig:TA}
\end{figure}

The simulated random field \eqref{clt1} is the sum of $L$ basic random fields (Gegenbauer waves), each of which varies along the meridians passing through a vector (pole) uniformly distributed on the sphere while it remains constant along the parallels orthogonal to this pole. We refer this construction as the `turning arcs' algorithm, by analogy with the turning bands method in which a random field in the Euclidean space is obtained by spreading basic random fields that varies along a direction spanned by a random vector and are constant along the hyperplanes orthogonal to this vector \citep{matheron1973intrinsic, mantoglou1982turning, lantu2002} (Figure \ref{fig:TA}).

As pointed out in \citet{emery2016} for the turning bands method, the process time of the turning arcs algorithm is, up to a pre-processing cost for generating the random vectors $\{\bm{\omega}_\ell: \ell = 1, \cdots, L\}$ and random variables $\{\kappa_\ell: \ell = 1, \cdots, L\}$, proportional to the number $L$ of basic random fields and to the number of target locations on the sphere, and turns out to be considerably fast. It is even faster than the spectral algorithms where the $L$ basic random fields are spherical harmonics or hyperspherical harmonics \citep{emery2019, lantuejoul2019}, insofar as the calculation of such harmonics is much more expensive than that of Gegenbauer polynomials, which can be easily computed by using \eqref{Gegenbauer}, see discussion in Section \ref{practicalaspects}.

\subsection{Extension to Vector-Valued Random Fields}

The goal of this section is to extend Proposition \ref{prop1} to the vector-valued case. Consider the sequence of Schoenberg matrices, $\{\bm{B}_{n,d}: n\in\mathbb{N}\}$, and the factorization $$\bm{B}_{n,d} = \bm{\Gamma}_{n,d} \, \bm{\Gamma}_{n,d}^\top, \qquad n\in\mathbb{N}.$$
For instance,  $\bm{\Gamma}_{n,d}$ can be the Cholesky factor of $\bm{B}_{n,d}$ or any square root of this matrix; in the latter case, $\bm{\Gamma}_{n,d}$ is symmetric since $\bm{B}_{n,d}$ is symmetric.
We use the notation $\bm{\gamma}^{(i)}_{n,d}$ for the $i$th column of $\bm{\Gamma}_{n,d}$. We observe that 
\begin{equation}
\label{gamma-matrix}
  \sum_{i=1}^p \bm{\gamma}^{(i)}_{n,d}\, [\bm{\gamma}^{(i)}_{n,d}]^\top    =  \bm{B}_{n,d}.
  \end{equation}
  
The following proposition provides a simulation algorithm for the vector-valued scenario. 
\begin{proposition}
\label{prop2}
 Let $\varepsilon$ be a random variable with zero mean and unit variance, $\bm{\omega}$ a random vector uniformly distributed on $\mathbb{S}^d$, $\iota$ a random integer uniformly distributed on $\{1, \hdots, p \}$ and $\kappa$ a random integer with $\mathbb{P}(\kappa = n) = a_n$, $n\in\mathbb{N}$, where $\{a_n: n\in\mathbb{N}\}$ is a probability mass sequence with a support containing that of the sequence of matrices $\{\bm{B}_{n,d}: n\in\mathbb{N}\}$. Suppose that all these random variables and vectors are independent. Then, 
\begin{enumerate}
\item[(1)] For $d=1$, the random field defined by
\begin{equation}
\label{simulacion-vector-d=1}
\bm{Z}(\bm{x}) = \varepsilon \, \sqrt{\frac{2p}{a_{\kappa}}} \, \bm{\gamma}^{(\iota)}_{\kappa,1} \,  \cos(\kappa \vartheta(\bm{\omega}, \bm{x})),        \qquad \bm{x}\in\mathbb{S}^1,
\end{equation}
is isotropic, with zero mean
%\begin{equation*}
%\label{media-vector-d=1}
%\mathbb{E}\{ \bm{Z}(\bm{x})  \} =  \sqrt{2 a_0} \sum_{i=1}^p \bm{\gamma}^{(i)}_{0,1}, \qquad \bm{x}\in\mathbb{S}^1,
%\end{equation*}
 and covariance function with isotropic part given by
%$$ \bm{K}(\vartheta)  =  2(1-a_0) \bm{B}_{0,1}  + \sum_{n=1}^\infty \bm{B}_{n,1} \cos(n \vartheta),    \qquad 0\leq \vartheta \leq \pi.$$
$$ \bm{K}(\vartheta)  =  \sum_{n=0}^\infty \bm{B}_{n,1} \cos(n \vartheta),    \qquad 0\leq \vartheta \leq \pi.$$

\item[(2)] For $d\geq 2$, the random field defined by
\begin{equation}
\label{simulacion-vector}
\bm{Z}(\bm{x}) = \varepsilon \, \sqrt{  \frac{ p (2\kappa + d -1)}{a_{\kappa} (d-1) } } \, \bm{\gamma}^{(\iota)}_{\kappa,d} \,  {G}_{\kappa}^{(d-1)/2}(\bm{\omega}^\top \bm{x}), \qquad \bm{x}\in\mathbb{S}^d,
\end{equation}
is isotropic, with zero mean
%\begin{equation*}
%\label{media-vector}
%\mathbb{E}\{ \bm{Z}(\bm{x})  \} = \sqrt{a_0} \sum_{i=1}^p \bm{\gamma}^{(i)}_{0,d}, \qquad \bm{x}\in\mathbb{S}^d,
%\end{equation*}
 and covariance function with isotropic part given by (\ref{schoenberg1_vector}).
%$$ \bm{K}(\vartheta) =  (1-a_0) \bm{B}_{0,d}  + \sum_{n=1}^\infty \bm{B}_{n,d} {G}_n^{(d-1)/2}(\cos \vartheta),    \qquad  0\leq \vartheta \leq \pi.$$
%$$ \bm{K}(\vartheta) =  \sum_{n=0}^\infty \bm{B}_{n,d} {G}_n^{(d-1)/2}(\cos \vartheta),    \qquad  0\leq \vartheta \leq \pi.$$
\end{enumerate} 

\end{proposition}

The proof of Proposition \ref{prop2} has been deferred to Appendix \ref{proof-prop2}. As for the scalar case, a central limit approximation of a vector-valued Gaussian random field is obtained by putting
\begin{equation}
\label{clt2}
\widetilde{\bm{Z}}(\bm{x}) =   \frac{1}{\sqrt{L}}\sum_{\ell = 1}^L {\bm{Z}}_\ell(\bm{x}), \qquad \bm{x}\in\mathbb{S}^d,
\end{equation}
where ${\bm{Z}}_1(\bm{x}),\hdots, {\bm{Z}}_L(\bm{x})$ are $L$ independent simulated copies, and $L$ is a large integer.

\subsection{Choice of the distributions of $\varepsilon$ and $\kappa$}
\label{distributionchoice}

The results presented in the previous subsections show that the desired spatial correlation structure is reproduced as soon as the random variable $\varepsilon$ has a zero mean and unit variance and the random integer $\kappa$ has a probability mass sequence whose support contains the support of the Schoenberg sequence associated with the covariance of the target random field. 

The choice of the distributions of $\varepsilon$ and $\kappa$ only impacts the rate of convergence of the central-limit approximation to the multivariate-Gaussian distribution. Which distributions yield a faster rate of convergence? To answer this question, following \cite{Chiles}, we focus on the marginal distribution of $\widetilde{Z}(\bm{x})$, as defined in \eqref{clt1} (the same exercise could be done in the multivariate case, by examining each component of $\widetilde{\bm{Z}}(\bm{x})$ as defined in \eqref{clt2}). The Berry-Ess\'{e}en inequality \citep{Berry, Esseen1} gives an upper bound for the Kolmogorov-Smirnov distance between the marginal distribution of $\widetilde{Z}(\bm{x})$ and a Gaussian distribution: 

\begin{equation}
    \label{berryesseen}
    \underset{z\in \mathbb{R}}{\sup} \enskip \Bigg\lvert \mathbb{P} \left( \frac{\widetilde{Z}(\bm{x})}{\sigma} < z \right) - G(z) \Bigg\rvert \leq \frac{\xi \mu_3^Z}{\sigma^3 \, \sqrt{L}},
\end{equation}

where $G$ is the standard Gaussian cumulative distribution function, $\mu_3^Z$ is the third-order absolute moment of the basic random field ${Z}(\bm{x})$ as defined in \eqref{simulacion_d=1} or \eqref{simulacion2}, that is: $\mu_3^Z = \mathbb{E}\{\lvert{Z}(\bm{x})\rvert ^3\}$, $L$ is the number of basic random fields as defined in \eqref{clt1}, $\sigma^2 = K(0)$ (variance of ${Z}(\bm{x})$ and $\widetilde{Z}(\bm{x})$) and $\xi$ is a constant greater than $0.4097$ and lower than $0.4748$ \citep{Esseen, Korolev, Shevtsova}.

Hereinafter, we focus on the case when $d \geq 2$ in order to express the third-order absolute moment $\mu_3^Z$ and to find out an upper bound for this moment. Accounting for the fact that $\varepsilon$ is independent of $\kappa$ and $\bm{\omega}$, one can write:
$$
\mu_3^Z = \mathbb{E}\left(\lvert \varepsilon \rvert ^3 \right) \mathbb{E} \Bigg\{ \left( \frac{b_{\kappa,d} (2\kappa+d-1)}{a_\kappa (d-1)} \right)^{3/2} \mathbb{E} \left(\lvert G_\kappa^{(d-1)/2}(\bm{\omega}^T \bm{x}) \rvert ^3 \, \Big | \kappa \right) \Bigg\}.
$$

For $\mu_3^Z$ to be minimum, the third-order absolute moment of $\varepsilon$ must be minimum. Jensen's moment inequality \citep{Jensen} implies that $\mathbb{E}\{ \lvert \varepsilon \rvert ^3 \} \geq \mathbb{E}\{ \varepsilon^2 \}^{3/2} = 1$, the equality being reached  when $\varepsilon$ has a symmetric two-point distribution concentrated at $-1$ and $+1$ (Rademacher distribution), i.e., $\varepsilon$ is a random sign with equal probability of being positive or negative. On the other hand, one has (Appendix \ref{proof-BerryEsseen}):
\begin{equation}
\label{mundG}    
\mu_3^G(n) := \mathbb{E} \left(\lvert G_\kappa^{(d-1)/2}(\bm{\omega}^T \bm{x}) \rvert ^3 \, \Big | \kappa = n \right) = \left\{
\begin{aligned}
& \mathcal{O}(n^{-3/2}) \text{ if $d=2$}\\
& \mathcal{O}(\ln n) \text{ if $d=3$}\\
& \mathcal{O}(n^{3\lfloor \frac{d-1}{2} \rfloor}) \text{ if $d\geq4$},
\end{aligned}
\right.
\end{equation}
where $\lfloor \cdot \rfloor$ denotes the floor function. Under these conditions, one has
\begin{equation}
\label{mu3}
\mu_3^Z = \frac{1}{(d-1)^{3/2}} \sum_{n} \frac{b_{n,d}^{3/2} \, (2n + d -1)^{3/2} \, \mu_3^G(n)}{a_n^{1/2}},
\end{equation}
the sum being extended over the integers $n$ such that $a_n>0$.

%and accounting for \eqref{Gegenbauer2}, one has (we here focus on the case when $d \geq 2$ and use \eqref{simulacion2}): 
%\begin{equation}
%\label{mu3}
%\mu_3^Z \leq \mathbb{E} \bigg\{ \left( \frac{b_{\kappa,d} (2\kappa + d -1)}{a_\kappa (d-1) } \right)^{3/2}  \left( \frac{\Gamma(d-1+\kappa)}{\Gamma(d-1) \Gamma(\kappa+1)} \right)^3 \bigg\} = \frac{S}{(d-1)^{3/2} \, \Gamma^{3}(d-1)} ,
%\end{equation}
%with
%\begin{equation}
%\label{S}
%S =  \sum_{n} \frac{b_{n,d}^{3/2} \, (2n + d -1)^{3/2} \, \Gamma^{3}(d-1+n)}{a_n^{1/2} \, \Gamma^{3}(n+1)},\\
%\end{equation}
%the sum being extended over the integers $n$ such that $a_n>0$. 

The following cases provide criteria to choose a probability mass sequence $\{a_n: n\in\mathbb{N}\}$ that yields a finite value for $\mu_3^Z$, therefore a finite upper bound in the Berry-Ess\'{e}en inequality \eqref{berryesseen}, ensuring the convergence of the distribution of $\widetilde{Z}(\bm{x})$ to a normal distribution with a rate in $L^{-1/2}$, where $L$ is defined in \eqref{clt1} or \eqref{clt2}:

\begin{enumerate}
    \item [Case 1.] The Schoenberg sequence $\{b_{n,d}: n\in\mathbb{N}\}$ has a finite support, i.e., $b_{n,d}$ is nonzero for finitely many values of $n$. In this case, any choice of the probability mass sequence $\{a_n: n\in\mathbb{N}\}$ leads to a finite value for $\mu_3^Z$, therefore to a finite upper bound in the Berry-Ess\'{e}en inequality.
    
    \item [Case 2.] The Schoenberg sequence $\{b_{n,d}: n\in\mathbb{N}\}$ is nonzero for infinitely many values of $n$ and is such that $\limsup_{n \to +\infty}{\sqrt[n]{b_{n,d}}} = r < 1$. In such a case, based on the Cauchy root convergence test, $\mu_3^Z$ is finite provided that the following condition holds:
    \begin{equation}
    \label{criterion2a}
    \liminf_{n \to +\infty} {\sqrt[n]{a_{n}}} \geq r^3.
    \end{equation}

  %  \item [Case 2.] The support of the Schoenberg sequence $\{b_{n,d}: n\in\mathbb{N}\}$ is $\mathbb{N}$ and the ratio of two consecutive Schoenberg coefficients tends to a value less than 1, i.e.
  %  \begin{equation}
  %  \label{criterion2a}
  %  \lim_{n \to +\infty} \frac{b_{n+1,d}}{b_{n,d}} = r < 1.
  %  \end{equation}
  %  In such a case, the ratio of two consecutive summands in \eqref{S} tends to a value less than 1, which implies the convergence of the sum \eqref{S} to a finite value, provided that the probability mass sequence $\{a_n: n\in\mathbb{N}\}$ fulfills the following condition:
  %  \begin{equation}
  %  \label{criterion2b}
  %  \lim_{n \to +\infty} \frac{a_{n+1}}{a_{n}}  = r^\prime \in ]r^{3},1[.
  %  \end{equation}

    \item [Case 3.] The Schoenberg sequence $\{b_{n,d}: n\in\mathbb{N}\}$ is nonzero for infinitely many values of $n$ and such that $b_{n,d} = \mathcal{O}(n^{-\theta})$. On the one hand, the convergence of the series $\{b_{n,d} \, G_n^{(d-1)/2}(1): n\in\mathbb{N}\}$ implies that $\theta$ is greater than $d-1$. On the other hand, using formula 6.1.46 of \citet{abramowitz}, it is found that the summand in \eqref{mu3} is $\mathcal{O}(a_n^{-1/2} \, n^{3/2-3\theta/2} \, \mu_3^G(n))$. Based on \eqref{mundG}, $\mu_3^Z$ is finite if $a_n \geq c \, n^{-\theta^\prime}$ when $n \geq n_0$, with $n_0 \in \mathbb{N}$, $c > 0$ and $\theta^\prime \in ]1,\theta_{\max}^\prime[$ with
    %If $a_n \geq c \, n^{-\theta^\prime}$ when $n \geq n_0$, with $n_0$ an integer and $c$ and $\theta^\prime$ positive scalars, then, based on \eqref{mundG}, $\mu_3^Z$ is finite if $1<\theta^\prime<\theta_{\max}^\prime$ with:
    \begin{equation}
    \label{thetaprimemax}
\theta_{\max}^\prime = \left\{
\begin{aligned}
& 3\theta-2 \text{ if $d=2$}\\
& 3\theta-5 \text{ if $d=3$}\\
& 3\theta-5-6 \Big\lfloor \frac{d-1}{2} \Big\rfloor \text{ if $d\geq4$}.
\end{aligned}
\right.
\end{equation}
    Because $\theta > d-1$, a value of $\theta^\prime$ can always be found in the nonempty interval $]1,\theta_{\max}^\prime[$ when $d=2$ and $d=3$. In contrast, for $d\geq 4$, $\theta$ must be greater than $2 \lfloor \frac{d+1}{2} \rfloor$ for the interval $]1,\theta_{\max}^\prime[$ to be nonempty.

\end{enumerate}

\section{Examples}
\label{examples}

\subsection{Example 1: Bivariate random field with negative binomial covariance on $\mathbb{S}^2$}

The isotropic negative binomial covariance with parameter $\delta \in ]0,1[$ and the associated Schoenberg sequence on the $2$-sphere are given by
\begin{equation}
\label{negativebinomial1}
K_{NB}(\vartheta;\delta) =  \frac{1-\delta}{\sqrt{1+\delta^2-2\delta \cos \vartheta}},     \qquad 0\leq \vartheta \leq \pi,
\end{equation}
\begin{equation}
\label{negativebinomial2}
b_{n,2}^{NB}(\delta) =  (1-\delta) \, \delta^n,     \qquad n \in \mathbb{N}.
\end{equation}

A bivariate negative binomial covariance model and its associated Schoenberg sequence can be obtained as follows \citep{emery2019}: 
\begin{equation}
\label{negativebinomial3}
\bm{K}_{NB}(\vartheta;\bm{\delta},\rho) = \left[
\begin{array}{ccc}
K_{NB}(\vartheta;\delta_{11}) & \rho\,K_{NB}(\vartheta;\delta_{12}) \\
\rho\,K_{NB}(\vartheta;\delta_{12}) & K_{NB}(\vartheta;\delta_{22})
\end{array}
\right], \qquad 0\leq \vartheta \leq \pi,
\end{equation}
\begin{equation}
\label{negativebinomial4}
\bm{B}_{n,2}^{NB}(\bm{\delta},\rho) = \left[
\begin{array}{ccc}
b_{n,2}^{NB}(\delta_{11}) & \rho\,b_{n,2}^{NB}(\delta_{12}) \\
\rho\,b_{n,2}^{NB}(\delta_{12}) & b_{n,2}^{NB}(\delta_{22})
\end{array}
\right],  \qquad n \in \mathbb{N},
\end{equation}
with $\bm{\delta} = (\delta_{11},\delta_{12},\delta_{22})$ such that $\delta_{11} < 1$, $\delta_{22} < 1$, $\delta_{12} \leq \min(\delta_{11},\delta_{22})$ and $\lvert \rho \rvert \leq \frac{\sqrt{(1-\delta_{11})(1-\delta_{22})}}{1-\delta_{12}}$.\\

\begin{figure}[htp]
    \begin{center}
    \includegraphics[width=5.5cm]{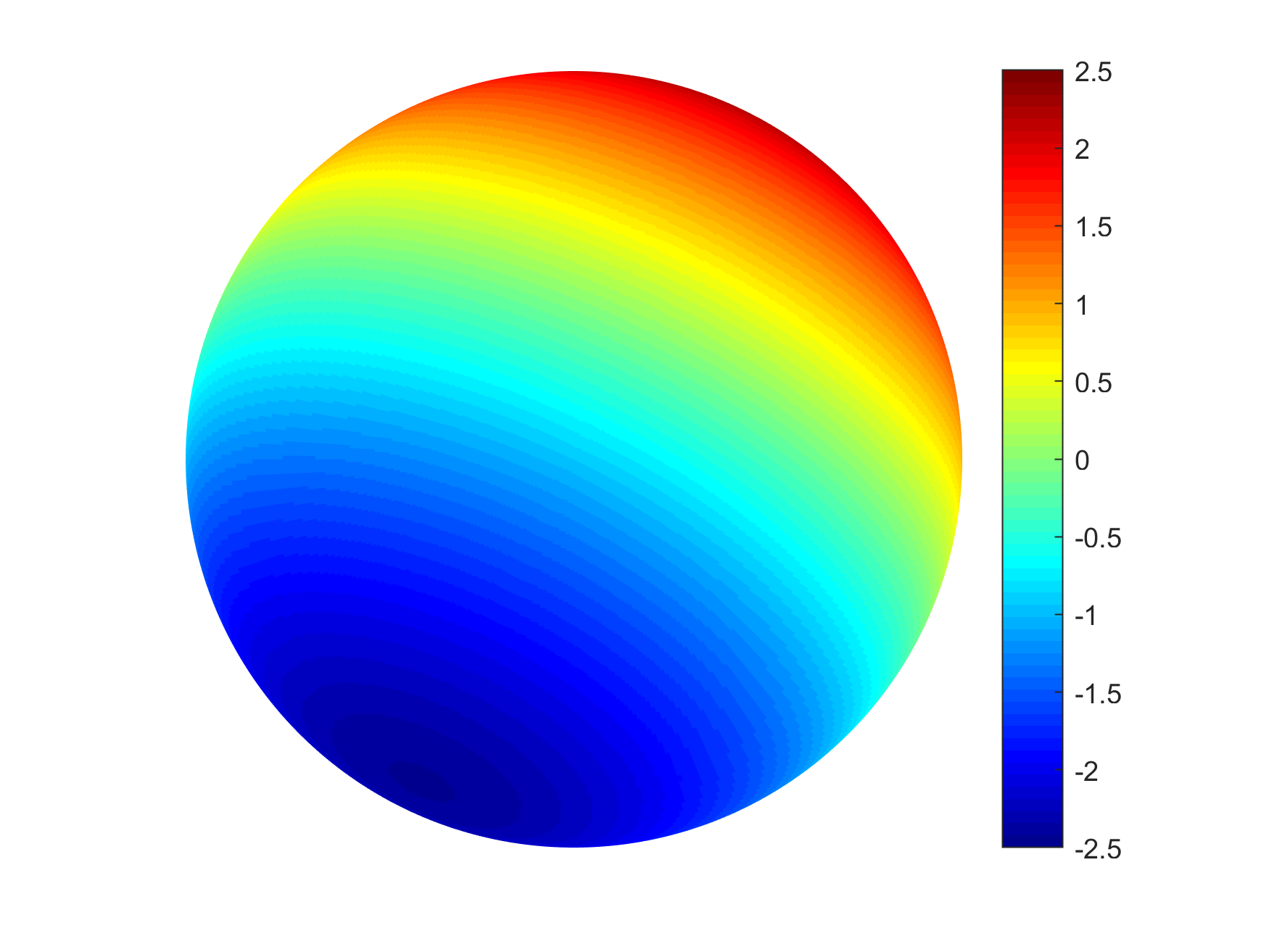} \hspace{0.1cm} \includegraphics[width=5.5cm]{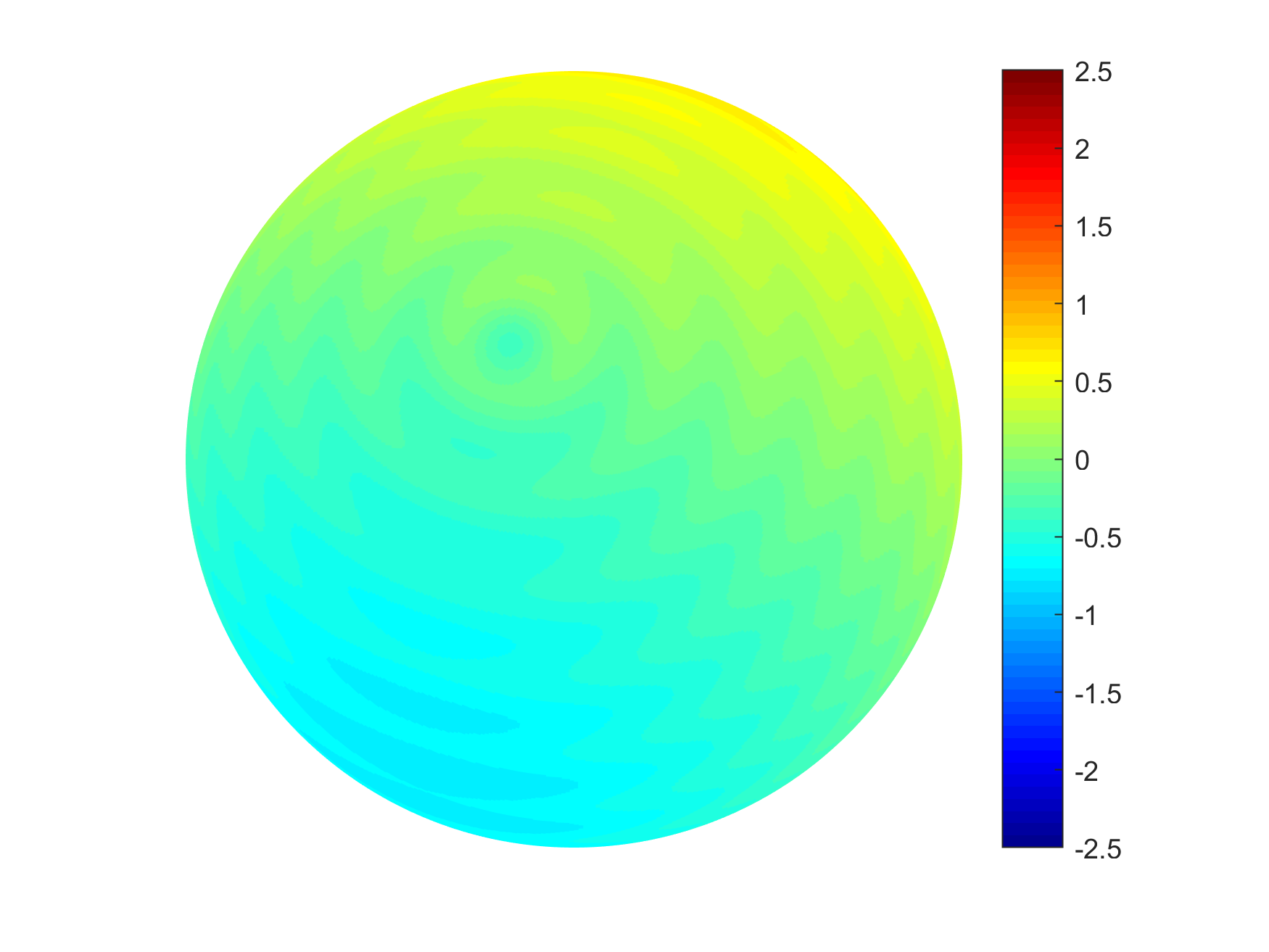} \hspace{0.1cm} \end{center}
    \begin{center}
    \includegraphics[width=5.5cm]{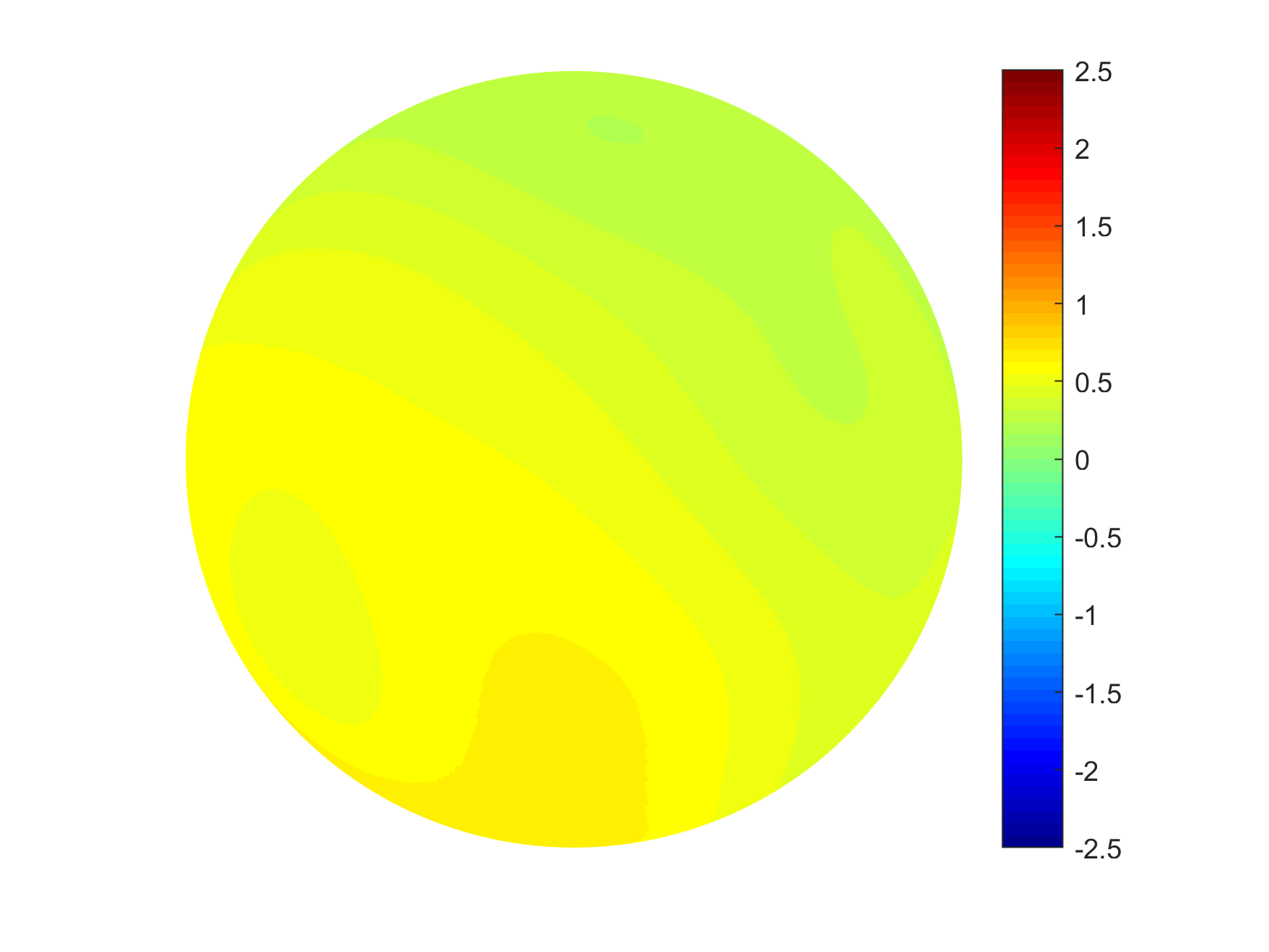} \hspace{0.1cm} \includegraphics[width=5.5cm]{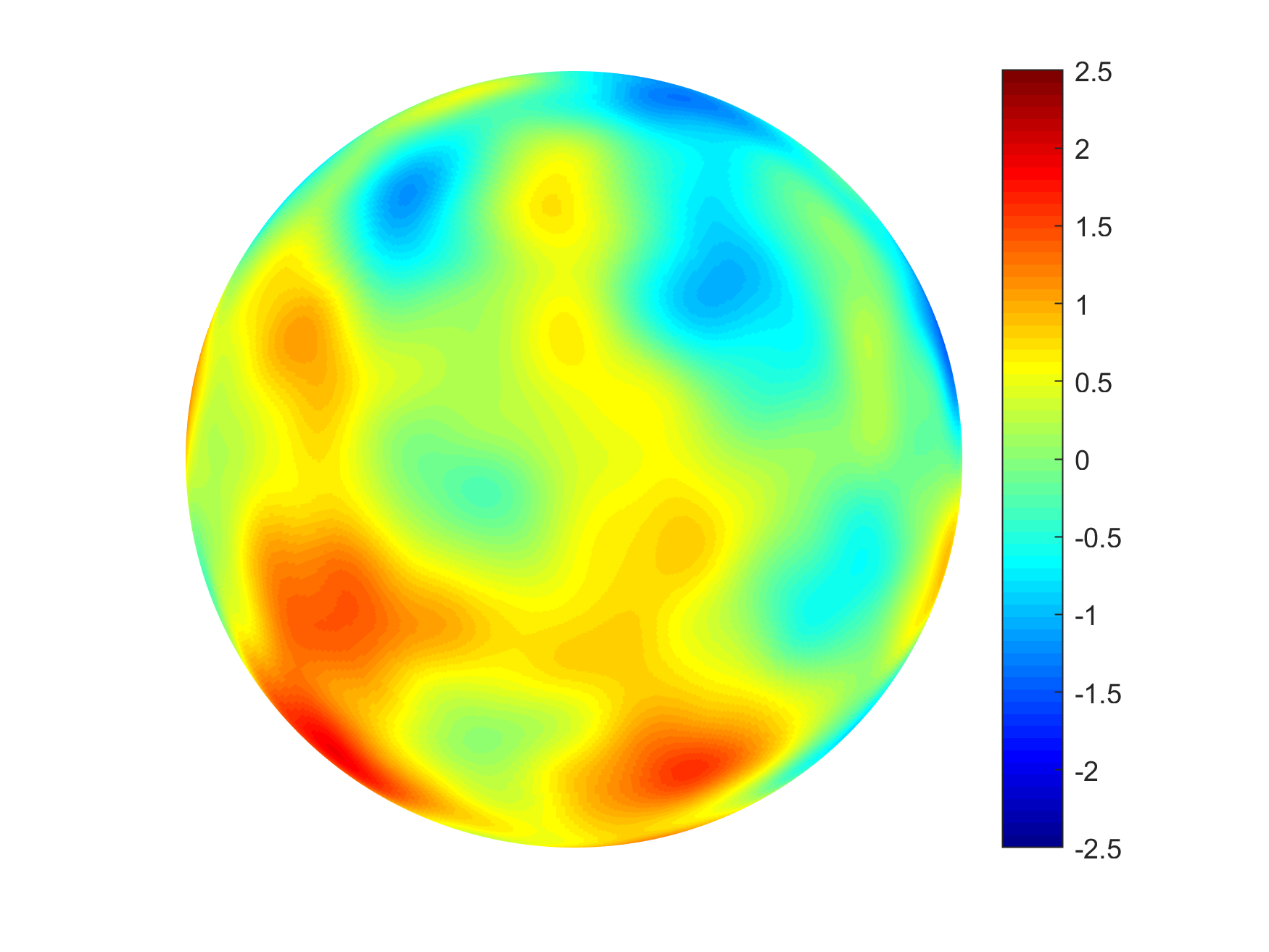} \hspace{0.1cm} \end{center}
    \begin{center}
    \includegraphics[width=5.5cm]{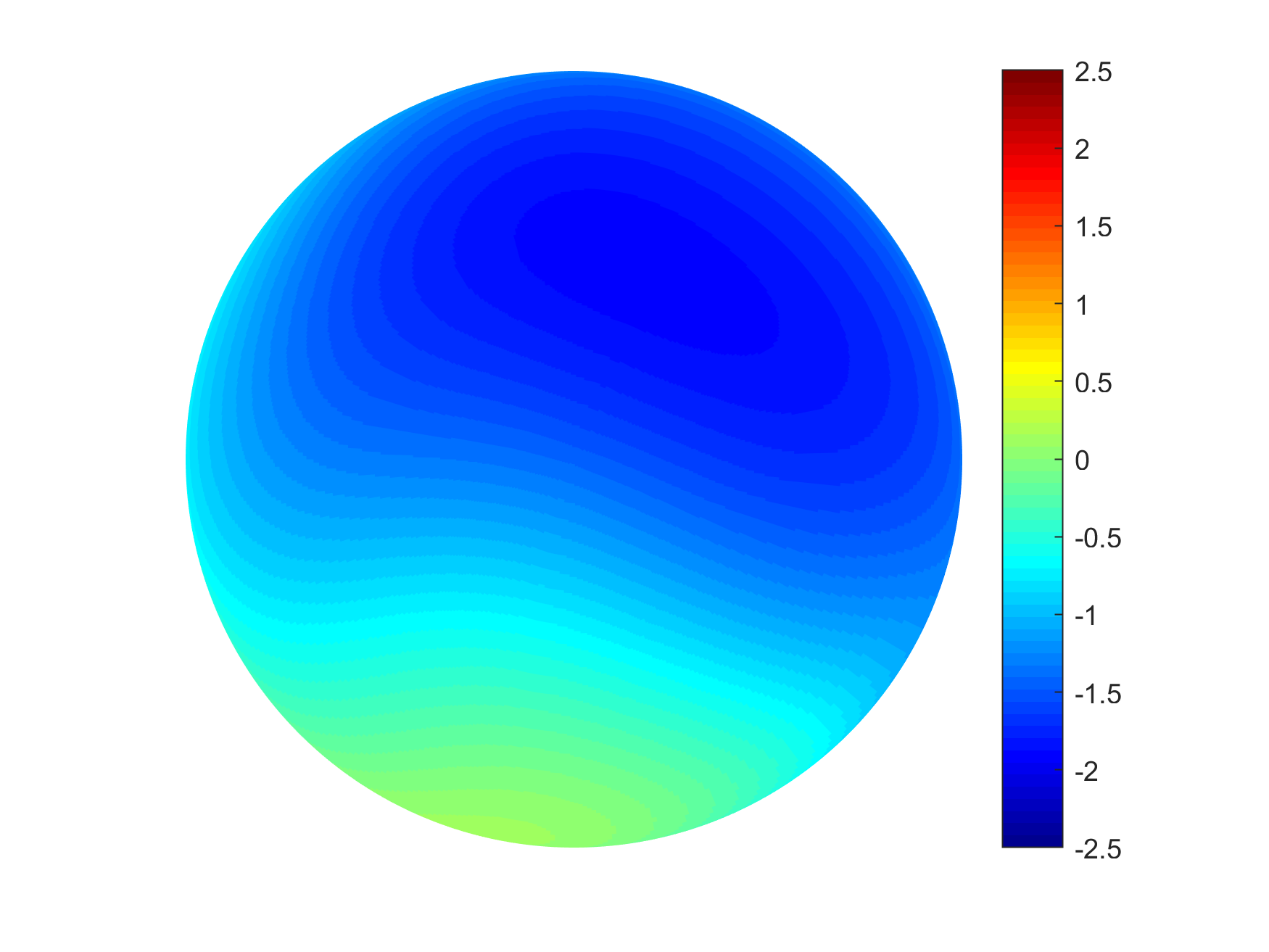} \hspace{0.1cm} \includegraphics[width=5.5cm]{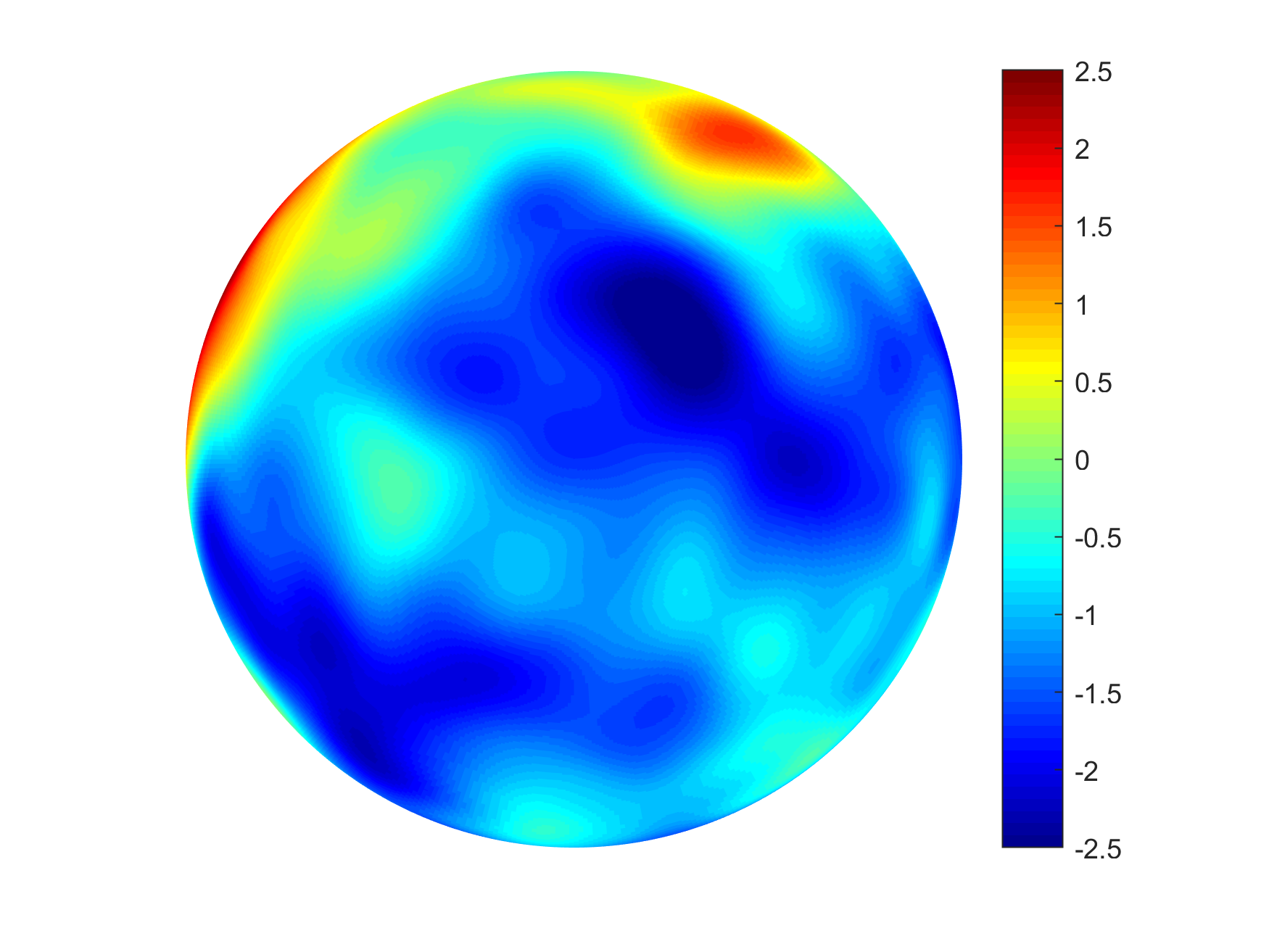} \hspace{0.1cm} \end{center}
\caption{Orthographic projections showing a realization of a bivariate random field with a negative binomial covariance ($\delta_{11}=\delta_{12}=0.2$, $\delta_{22}=0.7$ and $\rho=0.6$), obtained by using $L=15$ (top), $L=150$ (center) and $L=1500$ (bottom) basic random fields, a Rademacher distribution for $\varepsilon$ and a geometric distribution with success probability $0.01$ for $\kappa$. Left: first random field component; right: second random field component} 
\label{fig:BN}
\end{figure}

Since $\{b_{n,2}^{NB}(\delta): n\in\mathbb{N}\}$ in \eqref{negativebinomial2} is a geometric series, one has $\limsup_{n \to +\infty}{\sqrt[n]{b_{n,2}^{NB}(\delta)}} = \delta$. According to \eqref{criterion2a}, to ensure a finite Berry-Ess\'{e}en bound in \eqref{berryesseen} for both components of a bivariate random field with covariance \eqref{negativebinomial3}, it suffices to choose a probability mass sequence $\{a_{n}: n\in\mathbb{N}\}$ such that $\liminf_{n \to +\infty} {\sqrt[n]{a_{n}}} \geq \min(\delta_{11}^3,\delta_{22}^3)$. As an illustration, Figure \ref{fig:BN} shows orthographic projections of one realization of a bivariate random field obtained by applying the turning arcs algorithm with the following parameters:
\begin{itemize}
    \item $\delta_{11}=\delta_{12}=0.2$, $\delta_{22}=0.7$, $\rho=0.6$;
    \item $L=15$, $150$ or $1500$;
    \item $\varepsilon$ with a Rademacher distribution;
    \item $\kappa$ with a geometric distribution with success probability $0.01$;
    \item discretization of $\mathbb{S}^2$ into $500 \times 500$ faces with regularly-spaced colatitudes and longitudes.
\end{itemize}
Arc-shaped artifacts (striations) can be observed on the projections obtained with $L=15$, which indicates that the finite-dimensional distributions of the associated random field deviate from the multivariate-normal distributions expected for a Gaussian random field. This phenomenon is similar to the banding or striping effect of the continuous spectral and turning bands methods in the Euclidean space \citep{mantoglou1982turning, tompson, emery2006, emery2008}. The artifacts are no longer perceptible on the projections obtained with $L=150$ or $L=1500$ basic random fields, which display realizations that are visually close to that of a Gaussian random field, in agreement with the central limit theorem.  

\subsection{Example 2: Bivariate random field with spectral-Mat\'{e}rn covariance on $\mathbb{S}^2$}

The isotropic spectral-Mat\'{e}rn covariance with parameters $\alpha > 0$ and $\nu > 0$ on the $2$-sphere, hereafter denoted by $K_{SM}(\vartheta;\alpha,\nu)$ with $0 \leq \vartheta \leq \pi$, is associated with the following Schoenberg sequence \citep{guinness2016isotropic}:
\begin{equation}
\label{spectralmatern1}
b_{n,2}^{SM}(\alpha,\nu) =  \frac{(n^2+\alpha^2)^{-\nu-1/2}}{\sum_{k=0}^{+\infty}(k^2+\alpha^2)^{-\nu-1/2}},     \qquad n \in \mathbb{N}.
\end{equation}

As $n$ gets very large, the Schoenberg coefficient $b_{n,2}^{SM}(\alpha,\nu)$ is asymptotically of the order of $n^{-\theta}$ with $\theta=2\nu+1$. Based on the third case presented in Section \ref{distributionchoice}, a finite Berry-Ess\'{e}en bound is obtained when the probability mass sequence $\{a_{n}: n\in\mathbb{N}\}$ has a zeta distribution with parameter $\theta^{\prime} \in ]1,6\nu+1[$ (Eq. \eqref{thetaprimemax}), i.e.,   
\begin{equation}
\label{spectralmatern3}
a_n = \frac{1}{\zeta(\theta^{\prime})} n^{-\theta^{\prime}},
\end{equation}
where $\zeta$ refers to the Riemann zeta function \citep{abramowitz}. The simulation of a random variable $\kappa$ with such a zeta distribution can be done by the acceptance-rejection algorithm proposed by \cite{Devroye}.

A bivariate spectral-Mat\'{e}rn covariance model and its associated Schoenberg sequence can be obtained as follows \citep{emery2019}: 
\begin{equation}
\label{spectralmatern2}
\bm{K}_{SM}(\vartheta;\alpha,\bm{\nu},\rho) = \left[
\begin{array}{ccc}
K_{SM}(\vartheta;\alpha,\nu_{11}) & \rho\,K_{SM}(\vartheta;\alpha,\nu_{12}) \\
\rho\,K_{SM}(\vartheta;\alpha,\nu_{12}) & K_{SM}(\vartheta;\alpha,\nu_{22})
\end{array}
\right], \qquad 0\leq \vartheta \leq \pi,
\end{equation}
with $\alpha > 0$, $\bm{\nu} = (\nu_{11},\nu_{12},\nu_{22})$, $\nu_{11} > 0$, $\nu_{22} > 0$, $\nu_{12} \geq \frac{\nu_{11}+\nu_{22}}{2}$ and $\lvert \rho \rvert \leq \min(1,\alpha^{2\nu_{12}-\nu_{11}-\nu_{22}})$. \\

\begin{figure}[htp]
    \begin{center}
    \includegraphics[width=5.5cm]{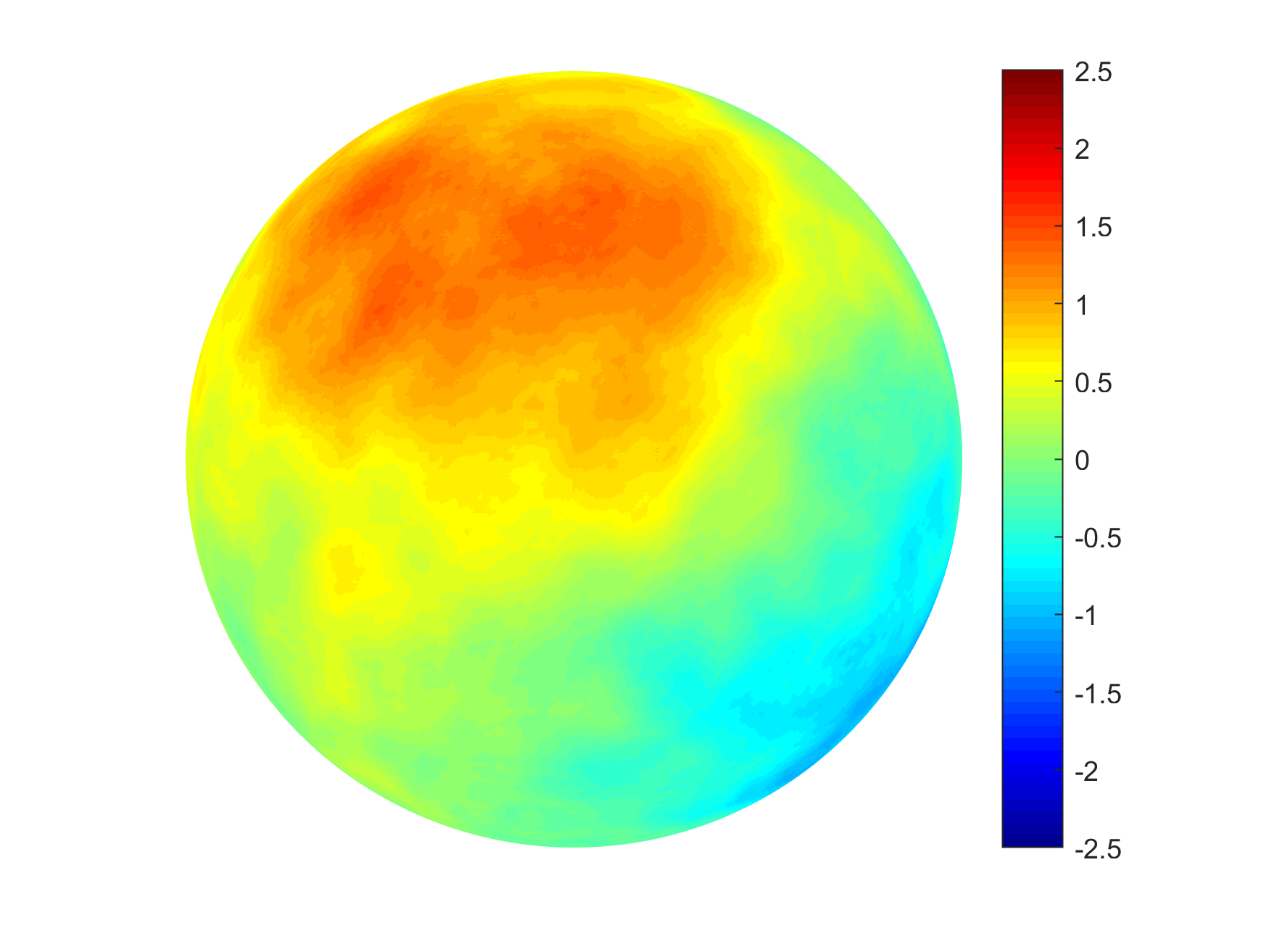} \hspace{0.1cm} \includegraphics[width=5.5cm]{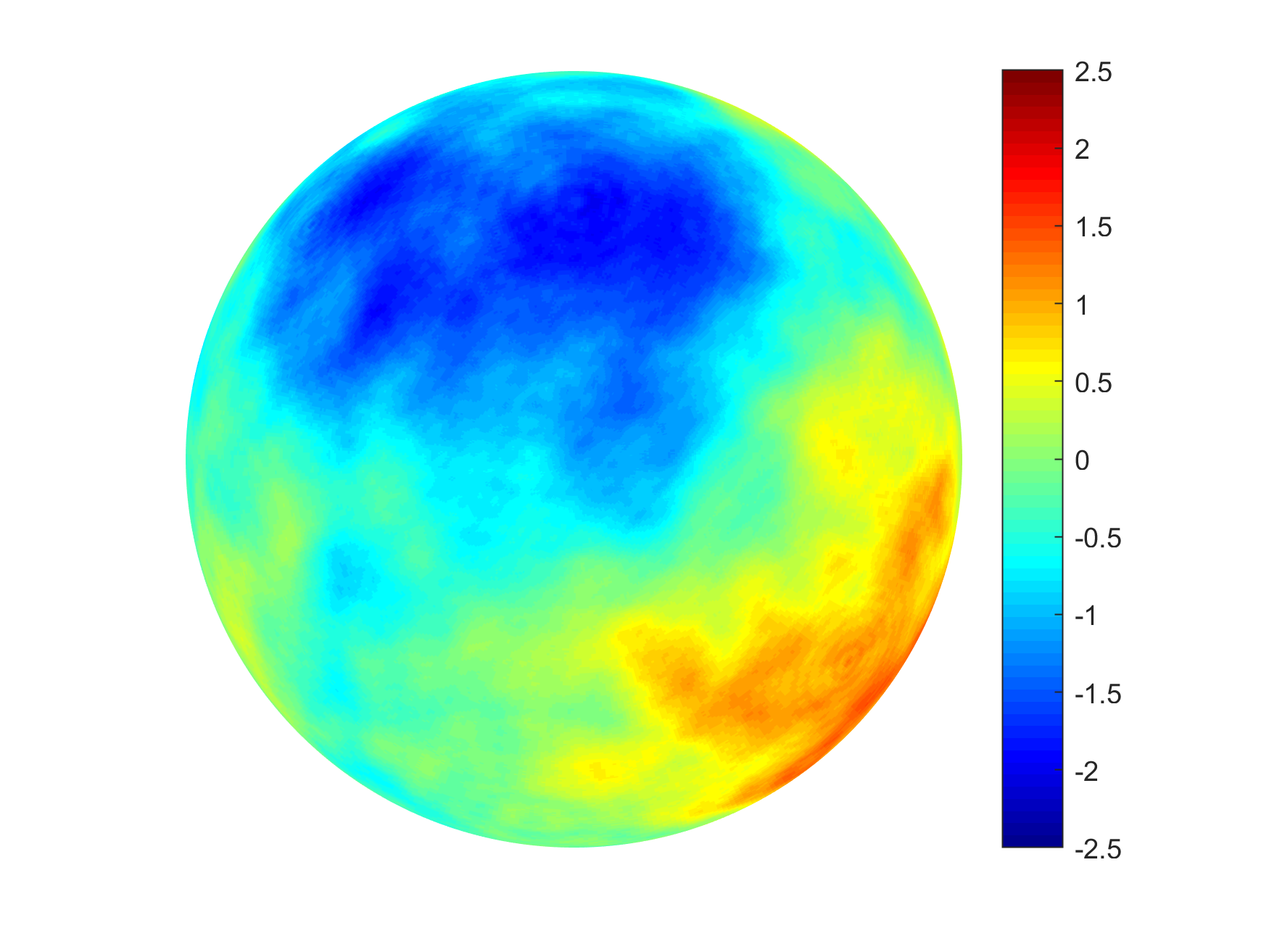} \hspace{0.1cm} \end{center}
\caption{Orthographic projections showing a realization of a bivariate random field with a spectral-Mat{\'e}rn covariance ($\alpha=1$, $\nu_{11}=2$, $\nu_{12}=\nu_{22}=0.75$ and $\rho=-0.6$), obtained by using $L=1500$ basic random fields, a Rademacher distribution for $\varepsilon$ and a zeta distribution with parameter $2$ for $\kappa$. Left: first random field component; right: second random field component} 
\label{fig:SM}
\end{figure}

The following illustration (Figure \ref{fig:SM}) shows orthographic projections of one realization of a bivariate random field obtained by applying the turning arcs algorithm with the following parameters:
\begin{itemize}
    \item $\alpha=1$, $\nu_{11}=2$, $\nu_{12}=\nu_{22}=0.75$, $\rho=-0.6$;
    \item $L=1500$;
    \item $\varepsilon$ with a Rademacher distribution;
    \item $\kappa$ with a zeta distribution with parameter $2$;
    \item discretization of $\mathbb{S}^2$ into $500 \times 500$ faces with regularly-spaced colatitudes and longitudes.
\end{itemize}

The two components are negatively correlated ($\rho<0$), the first one being smoother than the second one because $\nu_{11}>\nu_{22}$ \citep{guinness2016isotropic}. The striation effect is slightly perceptible in the right-hand side figure, which can be explained because the rate of convergence of the Schoenberg sequence $\{b_{n,2}^{SM}(\alpha,\nu_{22}): n\in\mathbb{N}\}$ is slower than that of the sequence $\{b_{n,2}^{SM}(\alpha,\nu_{11}): n\in\mathbb{N}\}$, hence the third-order absolute moment \eqref{mu3} and the upper bound in the Berry-Ess\'{e}en inequality \eqref{berryesseen} are higher: for the same number $L$ of basic random fields, the deviations from marginal normality and, \emph{a fortiori}, from multivariate normality, are likely to be more important for the second random field component than for the first one. 

%Note that, for a spectral-Mat{\'e}rn covariance with $\nu \leq 0.5$, the Berry-Ess\'{e}en upper bound \eqref{berryesseen} would be infinite irrespective of the choice of the distribution of $\kappa$ (the $a_n$'s obtained with \eqref{criterion0} do not constitute a convergent series), so that the rate of convergence to normality would remain unknown.

\subsection{Example 3: Univariate random field with generalized $\cal{F}$-covariance on $\mathbb{S}^3$}

%The isotropic $\cal{F}$ covariance has the following expression:
%\begin{equation}
%$\label{cal1}
%K_{\cal{F}}(\vartheta;\alpha,\nu,\tau) =  \sum_{n=0}^{+\infty} b_n^{\cal{F}}(\alpha,\nu,\tau) (\cos \vartheta)^n,     \qquad 0\leq \vartheta \leq \pi,
%\end{equation}
%with

The isotropic generalized $\cal{F}$-covariance on $\mathbb{S}^d$ is associated with the Schoenberg sequence $\{b_{n,d}^{\cal{F}}: n\in\mathbb{N}\}$ defined as follows \citep{alegria2018family}:

\begin{equation}
\label{cal2}
b_{n,d}^{\cal{F}}(\alpha,\nu,\tau) =  \frac{B(\alpha,\nu+\tau)}{B(\alpha,\nu)} \, \frac{(\alpha)_n \, (\tau)_n}{(\alpha+\nu+\tau)_n \, n!},     \qquad n \in \mathbb{N},
\end{equation}

where $\alpha>0$, $\nu>0$, $\tau>0$, $B(\cdot,\cdot)$ is the beta function and $(a)_n$ denotes the Pochhammer symbol \citep{abramowitz}. 
%It can be shown \citep{alegria2018family} that $K_{\cal{F}}$ is a valid covariance model on $\mathbb{S}^d$ for any $d \in \mathbb{N}^{*}$. 

As $n$ increases, the Schoenberg coefficient $b_{n,d}^{\cal{F}}(\alpha,\nu,\tau)$ is of the order of $n^{-\nu-1}$. As for the previous example, this suggests the use of a probability mass sequence $\{a_{n}: n\in\mathbb{N}\}$ with a zeta distribution with parameter $\theta^{\prime} \in ]1,3\nu-2[$ (Eq. \eqref{thetaprimemax}). 

\begin{figure}[htp]
    \begin{center}
    \includegraphics[width=5.5cm]{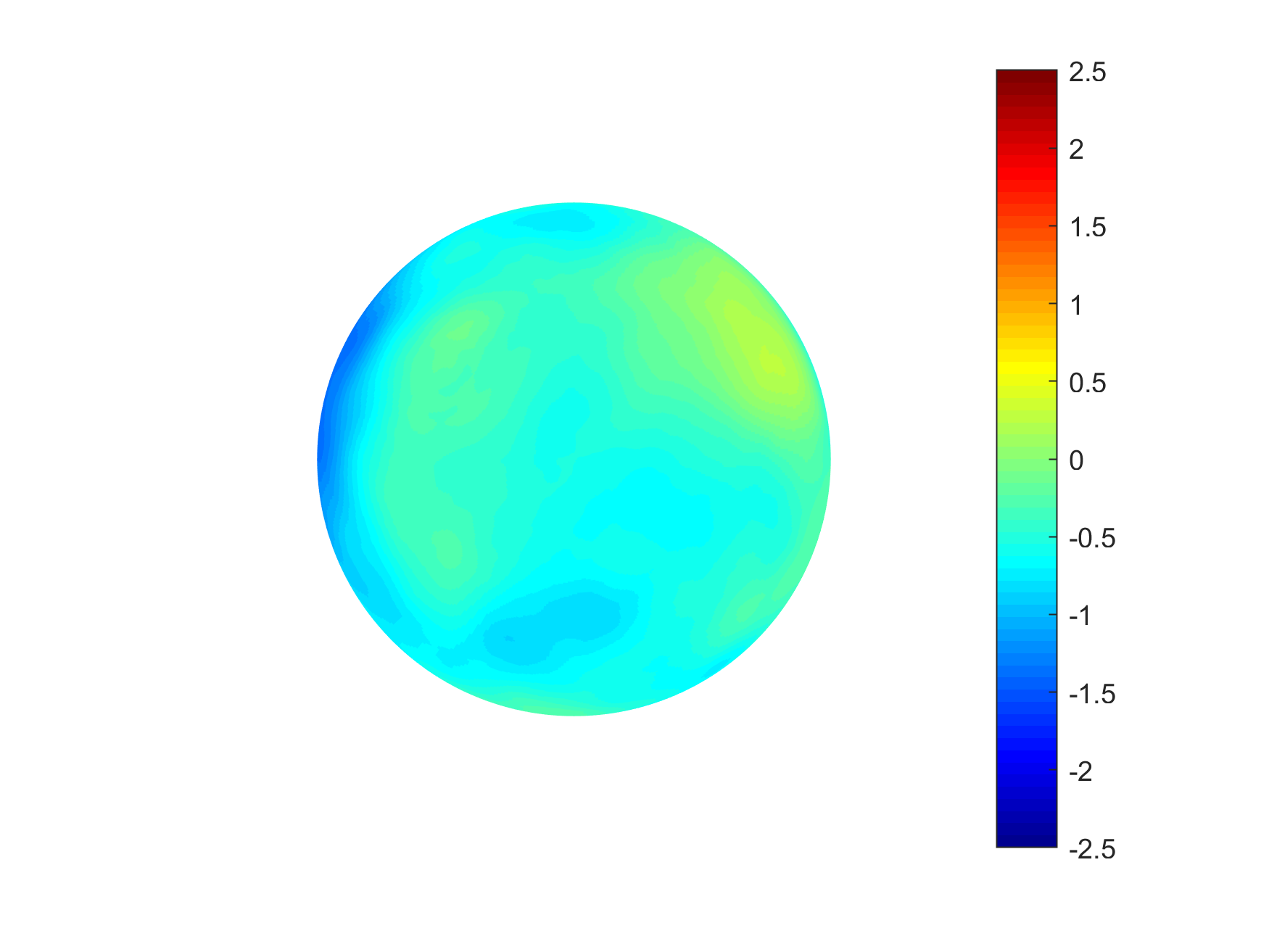} \hspace{0.1cm} \includegraphics[width=5.5cm]{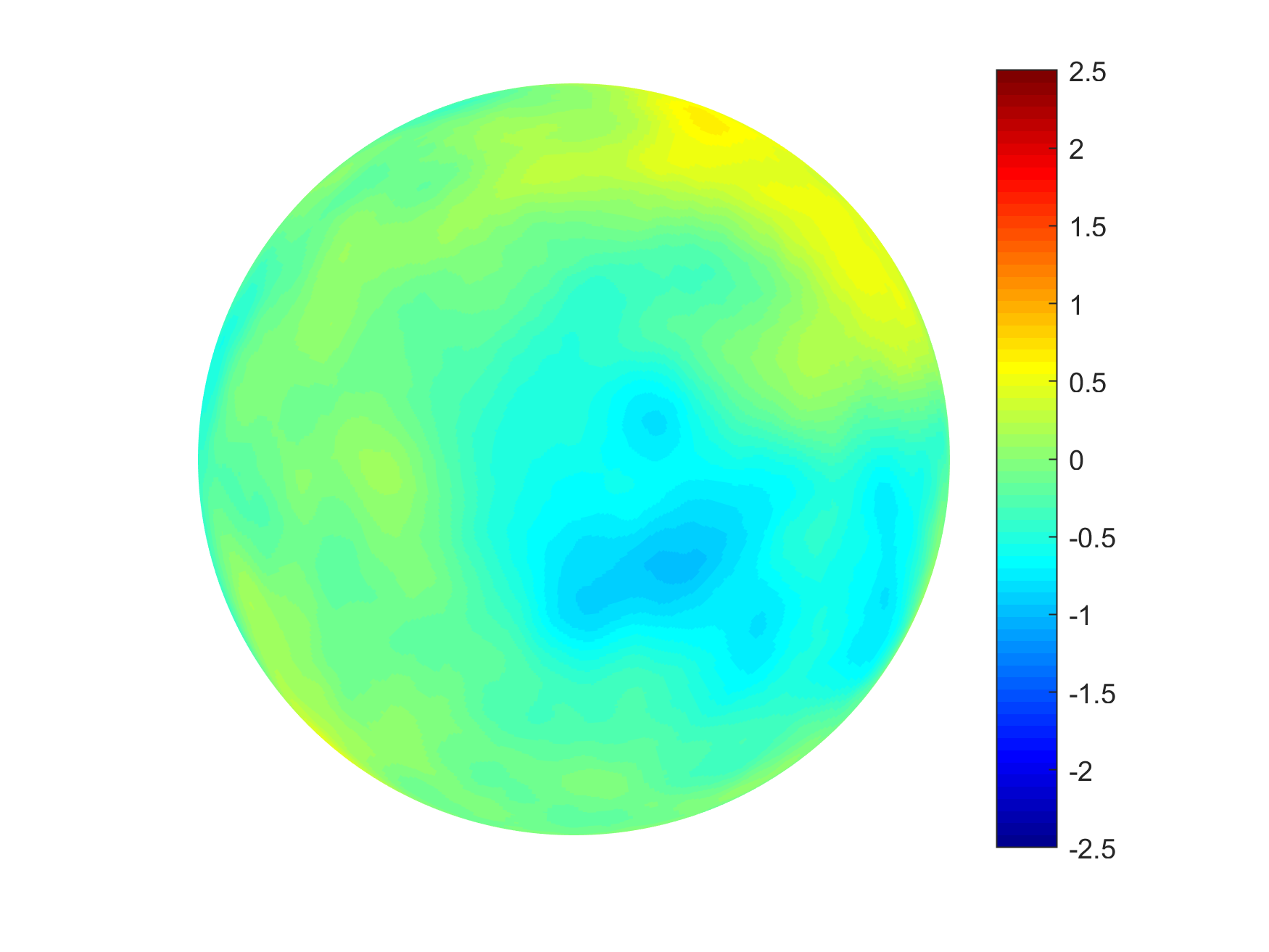} \hspace{0.1cm} \end{center}
    \begin{center}
    \includegraphics[width=5.5cm]{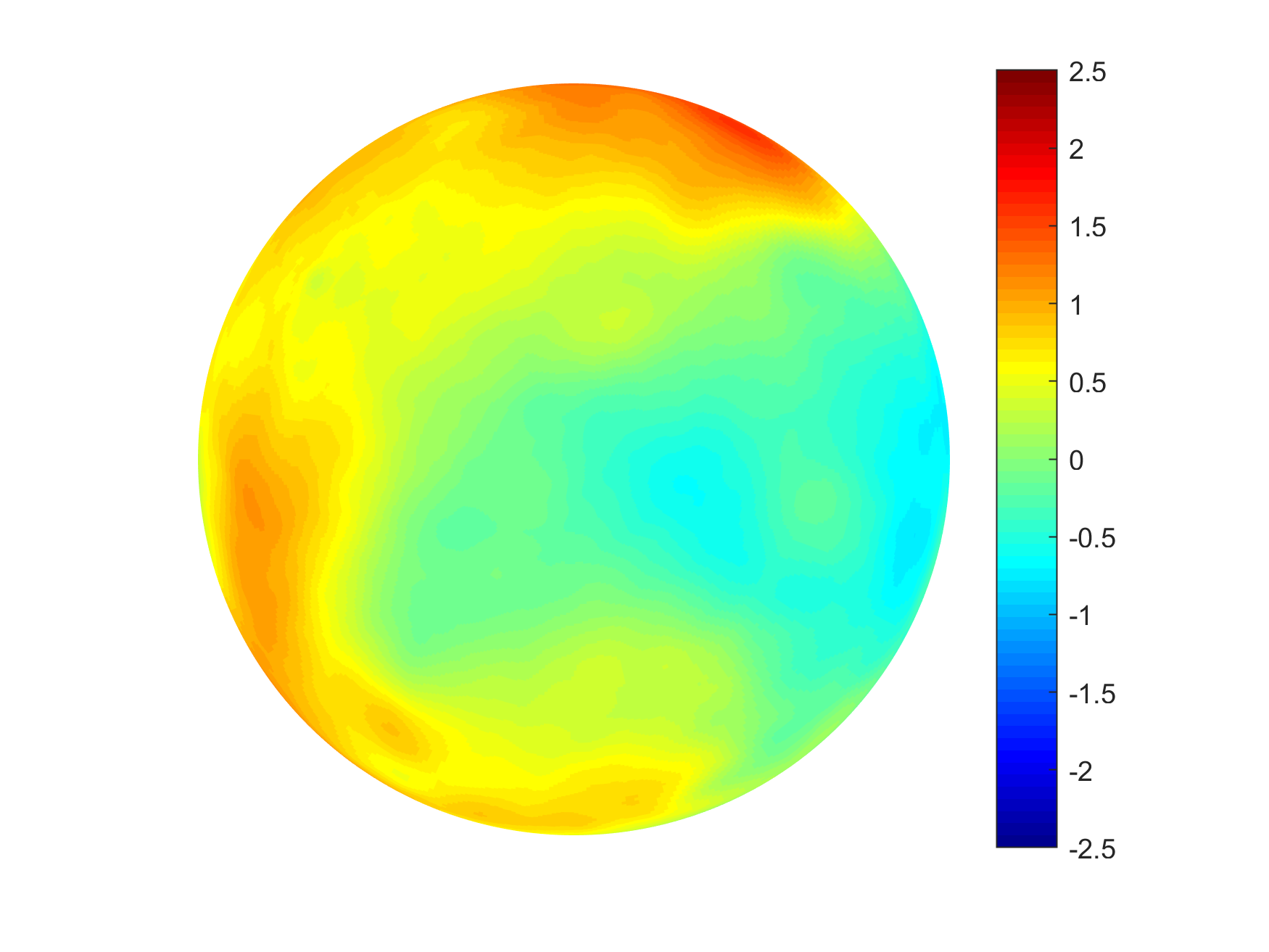} \hspace{0.1cm} \includegraphics[width=5.5cm]{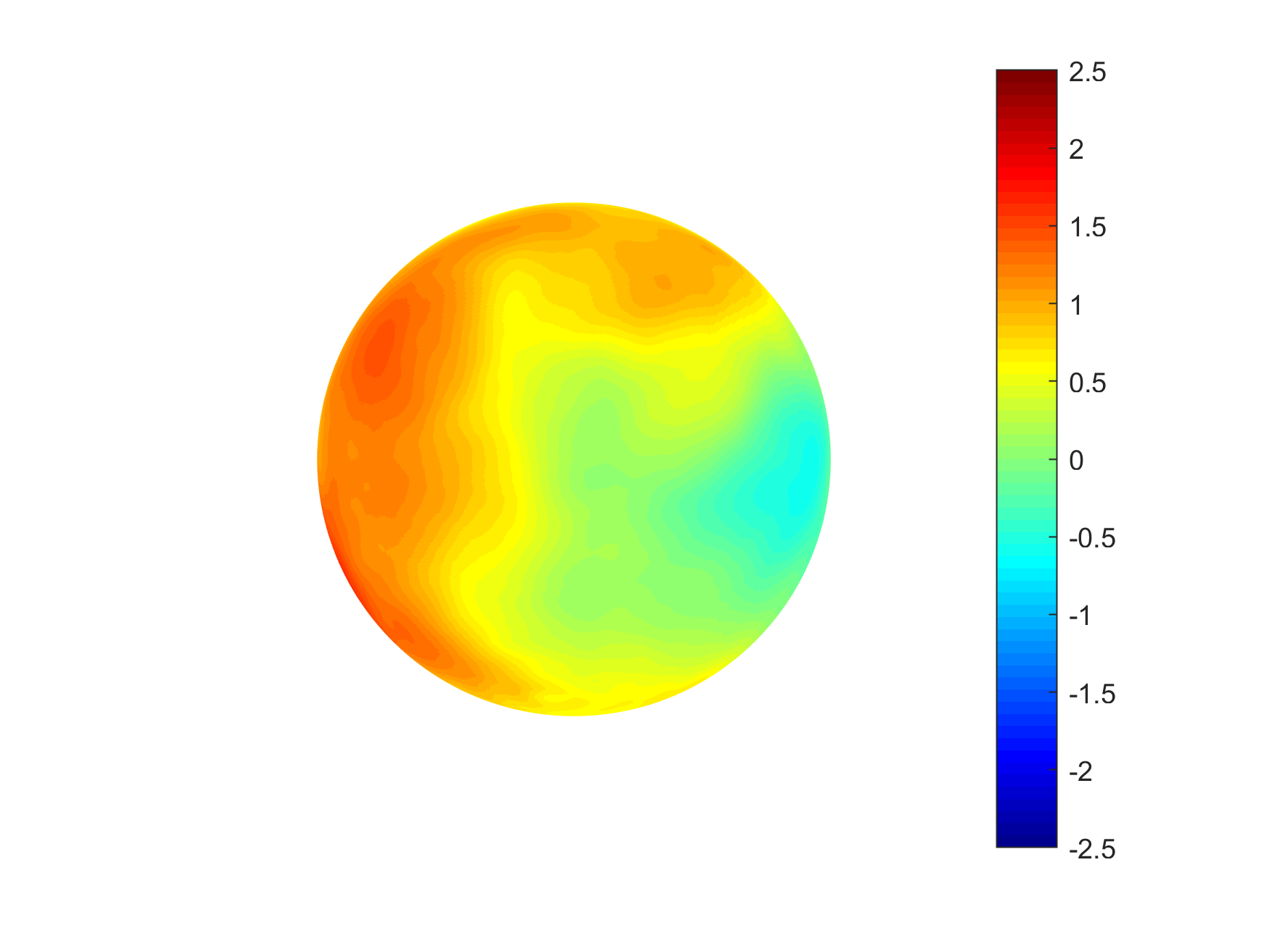} \hspace{0.1cm} \end{center}
\caption{Orthographic projections showing a realization
of a univariate random field with a  generalized $\cal{F}$-covariance ($\alpha=1$, $\nu=3.5$ and $\tau=2$) on the $3$-sphere, obtained by using $L=1500$ basic random fields, a Rademacher distribution for $\varepsilon$ and a zeta distribution for $\kappa$ with parameter $2$. Representations of the $2$-sphere corresponding to the sections of the $3$-sphere with the fourth coordinate equal to $-0.75$ (top left), $-0.25$ (top right), $0.25$ (bottom left) and $0.75$ (bottom right)
} 
\label{fig:Hypersphere}
\end{figure}

The following illustration (Figure \ref{fig:Hypersphere}) shows orthographic projections of one realization of a univariate random field obtained by applying the turning arcs algorithm with the following parameters:
\begin{itemize}
    \item $\alpha=1$, $\nu=3.5$, $\tau=2$;
    \item $L=1500$;
    \item $\varepsilon$ with a Rademacher distribution;
    \item $\kappa$ with a zeta distribution with parameter $2$;
    \item $d=3$;
    \item discretization of each $2$-sphere resulting from a cross-section of $\mathbb{S}^3$ into $500 \times 500$ faces with regularly-spaced colatitudes and longitudes.
\end{itemize}

%\begin{figure}[htp]
%    \begin{center}
%    \includegraphics[width=8cm]{F1.png} \hspace{0.1cm} \includegraphics[width=8cm]{F2.png} \hspace{0.1cm} \end{center}
%\caption{Orthographic projections showing a realization
%of a bivariate random field with $\cal{F}$ covariance ($\alpha_{11}=\alpha_{12}=1$, $\alpha_{22}=2$, $\nu_{11}=\nu_{12}=2$, $\nu_{22}=1$, $\tau_{11}=\tau_{12}=\tau_{22}=1$ and $\rho=0.7$), obtained by using $L=5000$  basic random fields, a standard normal distribution for $\varepsilon$ and a geometric distribution for $\kappa$ with success probability $0.01$} 
%\label{fig:F}
%\end{figure}

%\begin{figure}[htp]
    %\begin{center}
    %\includegraphics[width=8cm]{NB1.png} \hspace{0.1cm} \includegraphics[width=8cm]{NB2.png} \hspace{0.1cm} \end{center}
%\caption{Orthographic projections showing a realization
%of a bivariate random field with negative binomial covariance ($\delta_{11}=\delta_{12}=0.5$, $\delta_{22}=0.7$, $\tau_{11}=0.7$, $\tau_{12}=\tau_{22}=0.5$ and $\rho=0.8$), obtained by using $L=5000$  basic random fields, a standard normal distribution for $\varepsilon$ and a geometric distribution for $\kappa$ with success probability $0.01$} 
%\label{fig:NB}
%\end{figure}

\subsection{Example 4: Univariate random field with Chentsov covariance on $\mathbb{S}^d$}

\begin{figure}[htp]
    \begin{center}
    \includegraphics[width=5.5cm]{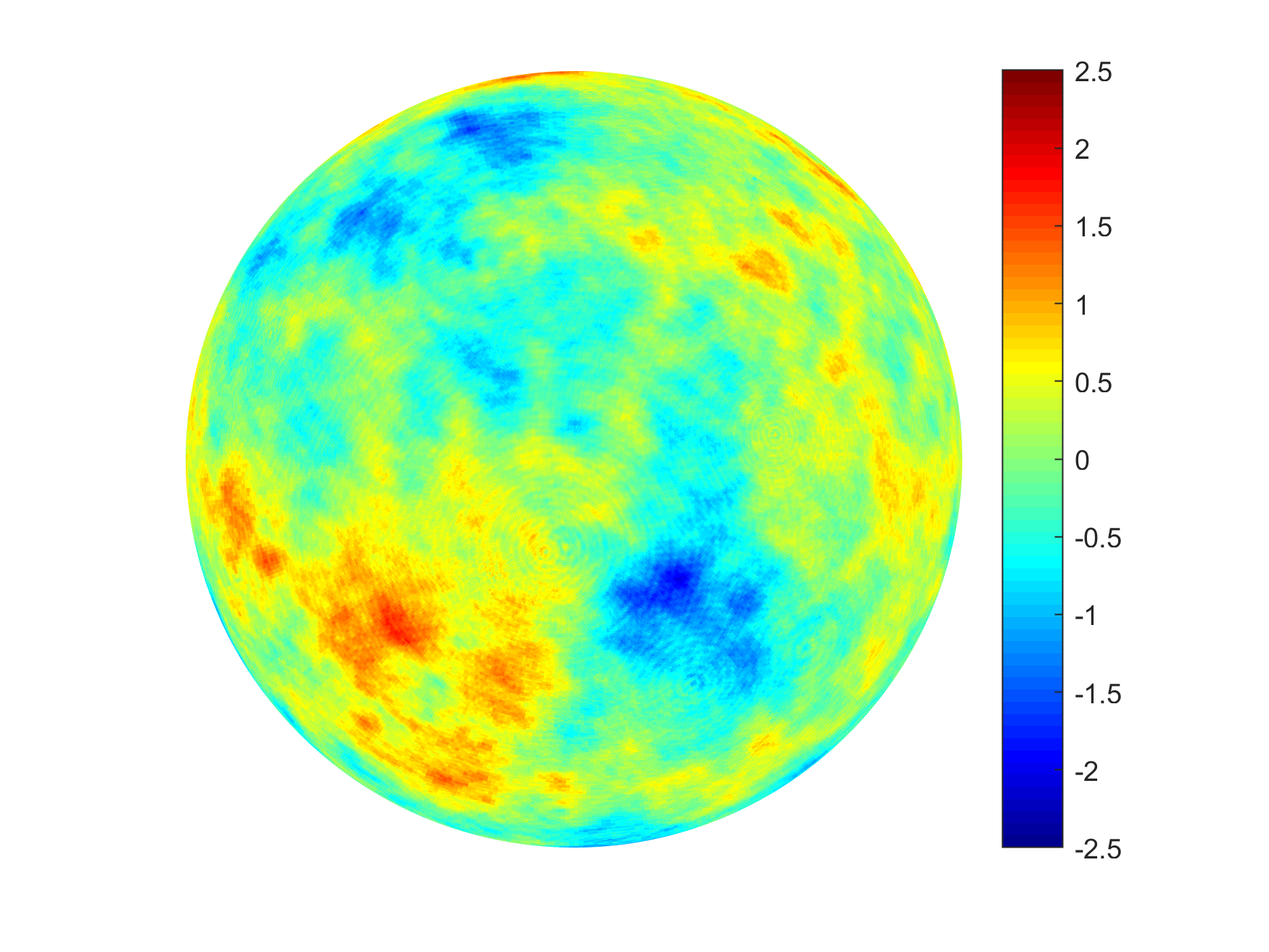} \hspace{0.1cm} \includegraphics[width=5.5cm]{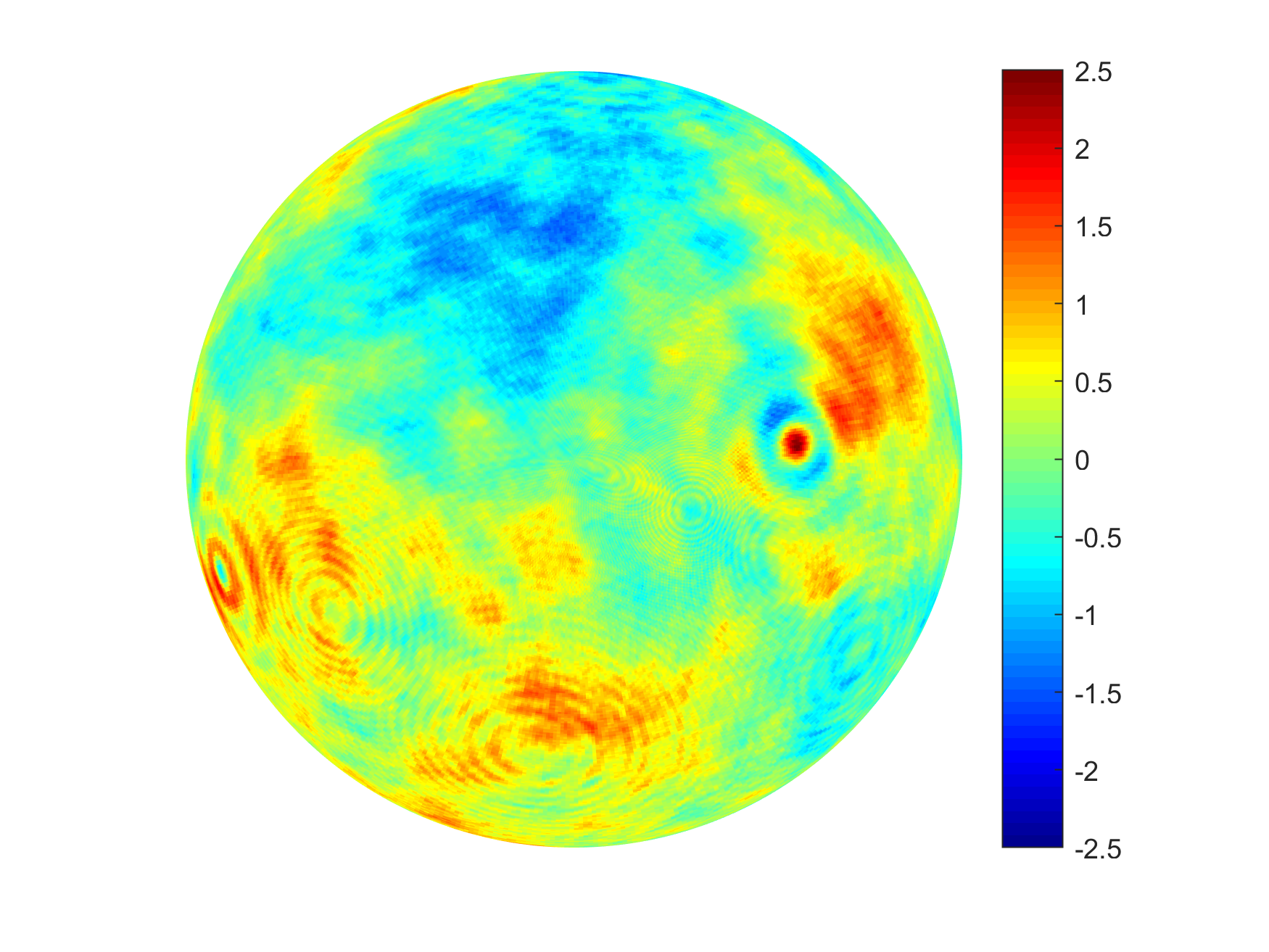} \hspace{0.1cm} \end{center}
    \begin{center}
    \includegraphics[width=5.5cm]{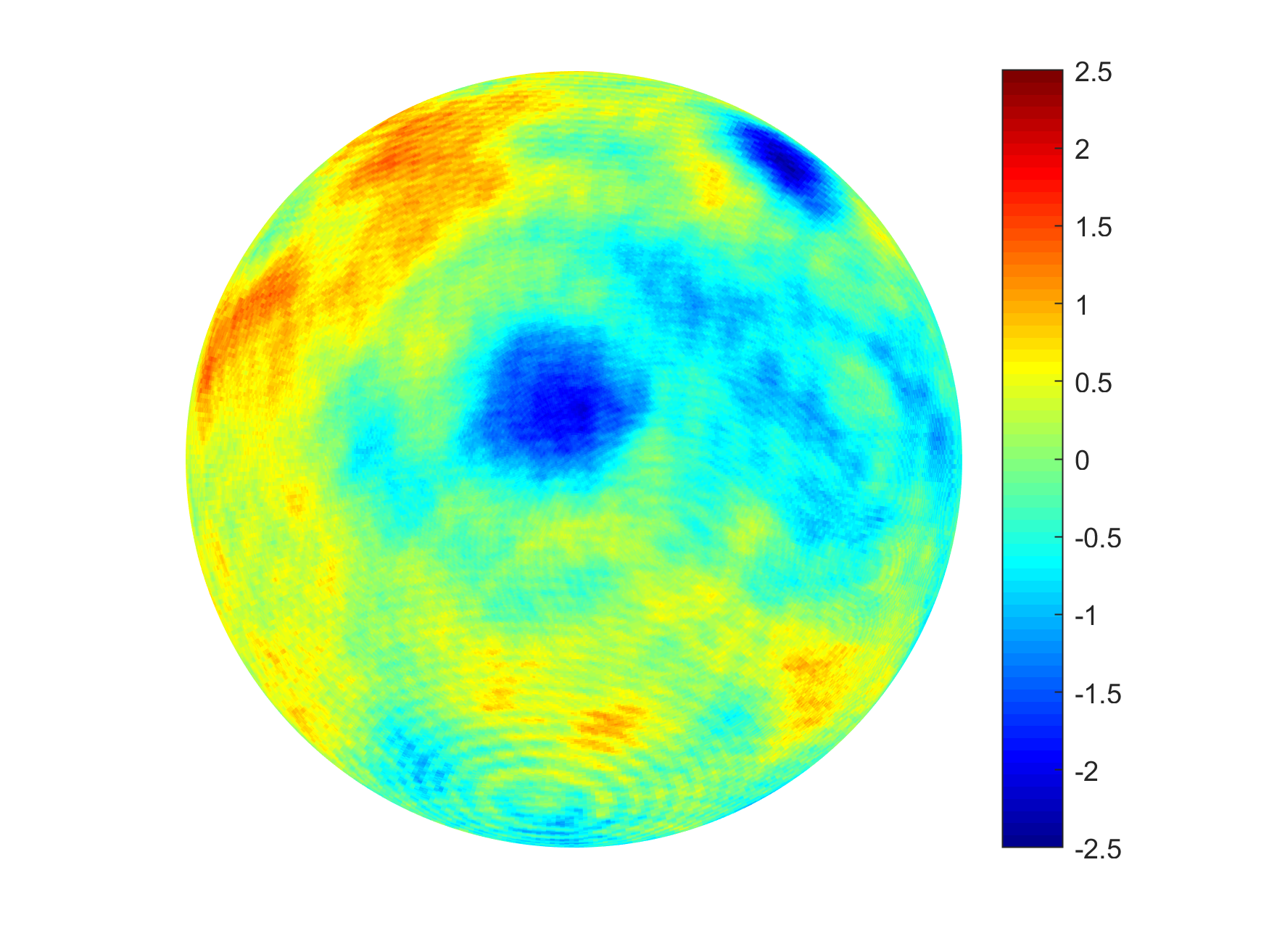} \hspace{0.1cm} \includegraphics[width=5.5cm]{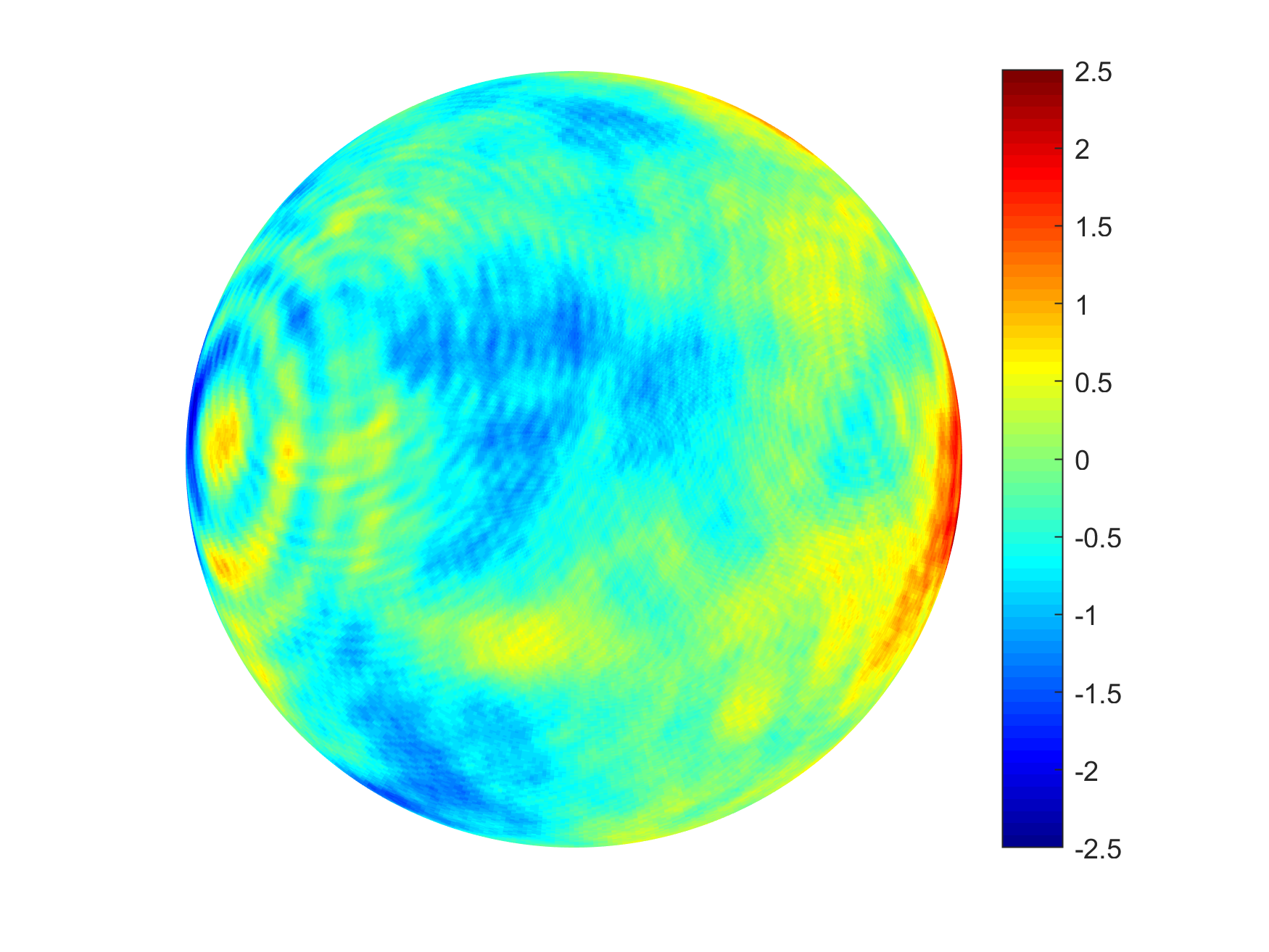} \hspace{0.1cm} \end{center}
    \begin{center}
    \includegraphics[width=5.5cm]{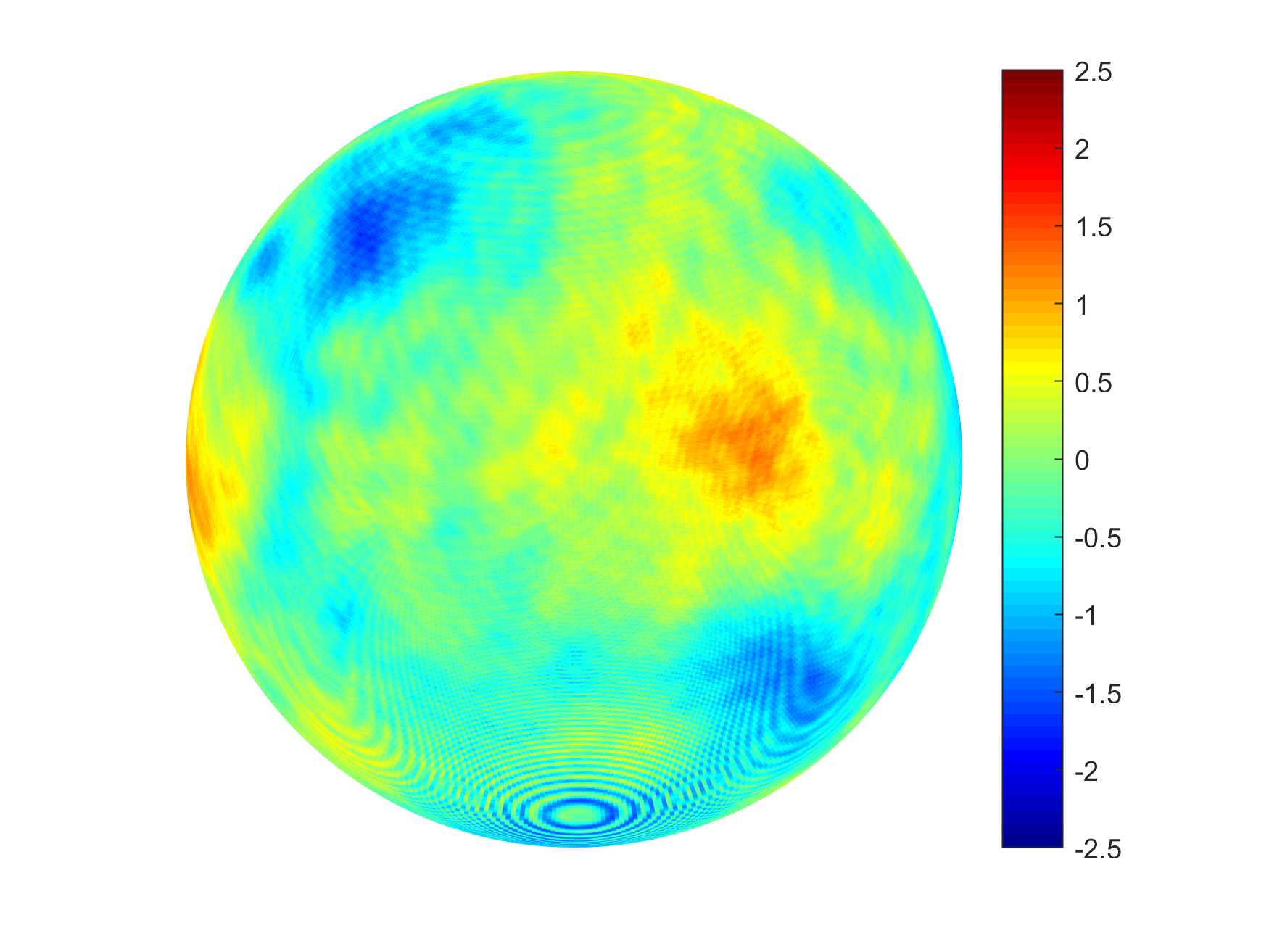} \hspace{0.1cm} \includegraphics[width=5.5cm]{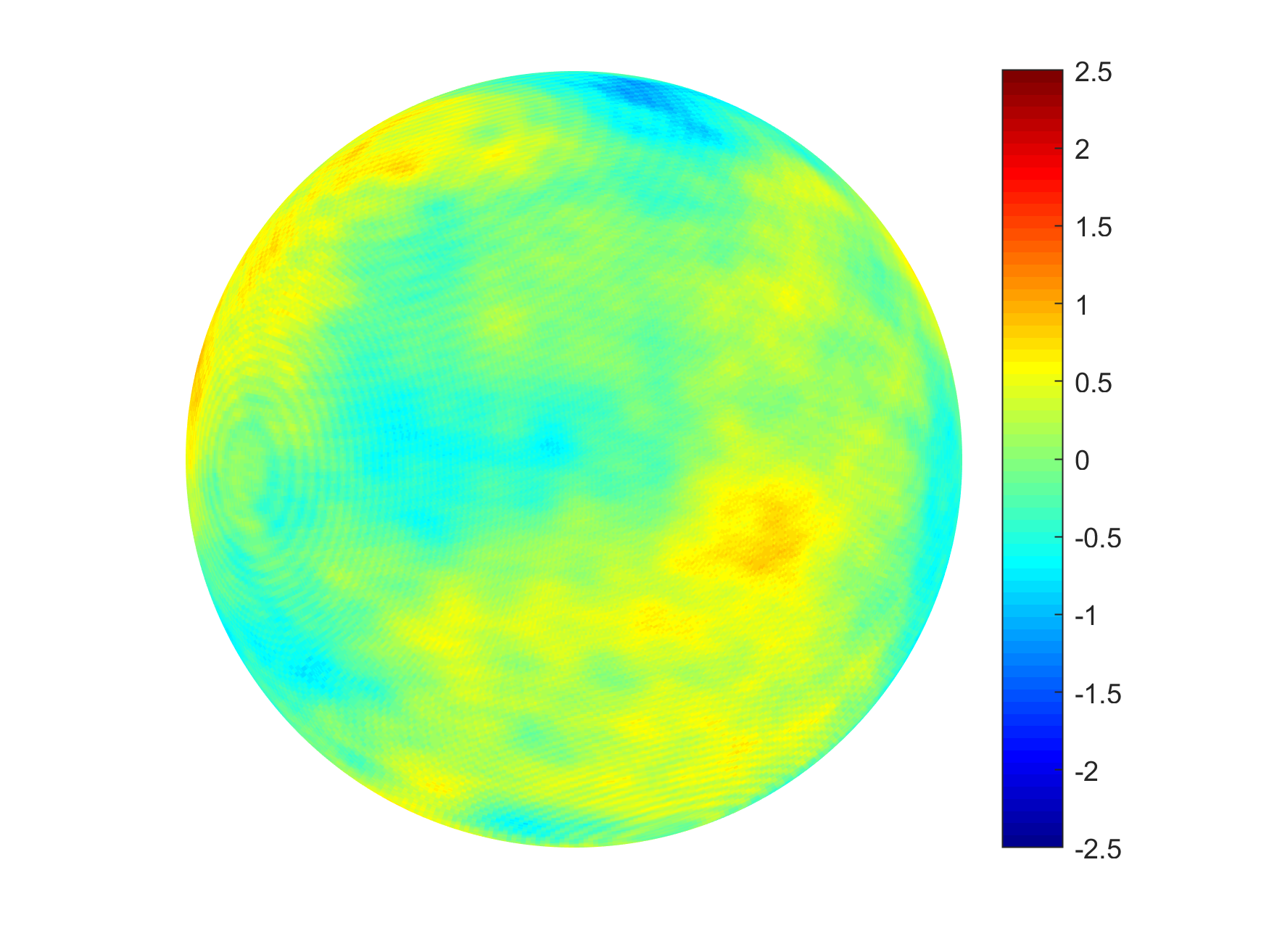} \hspace{0.1cm} \end{center}
    \begin{center}
    \includegraphics[width=5.5cm]{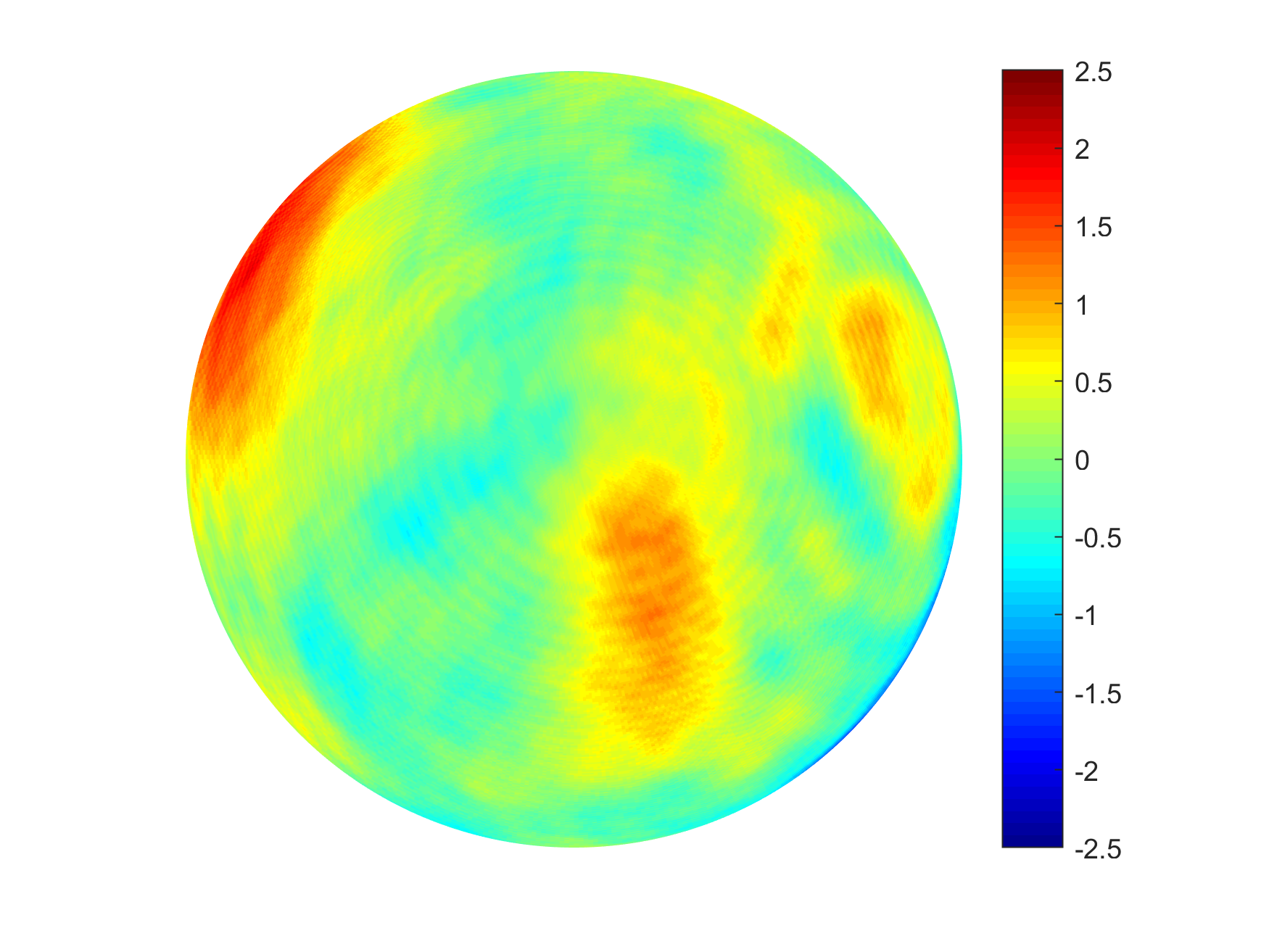} \hspace{0.1cm} \includegraphics[width=5.5cm]{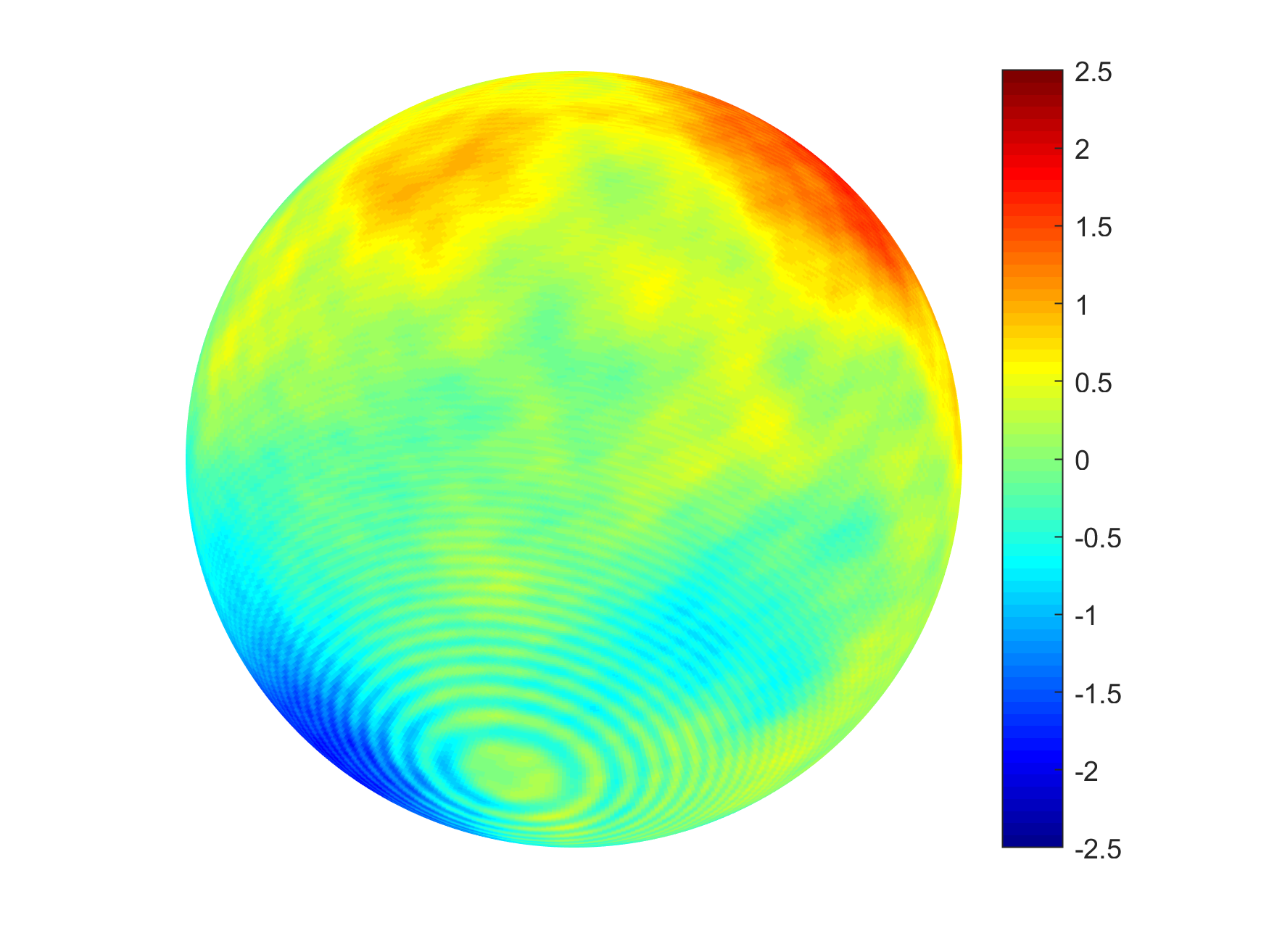} \hspace{0.1cm} \end{center}
    \caption{Orthographic projections showing a realization of a univariate random field with a Chentsov covariance, obtained by using $L=1500$ basic random fields, $\varepsilon$ with a Rademacher distribution and $\kappa$ of the form $2n+1$ with a zeta distribution of parameter $2$ for $n$. Representations of the $2$-sphere corresponding to the sections of the $d$-sphere with the last $d-2$ coordinates equal to 0. From top to bottom and left to right: $d=2$, $4$, $8$, $16$, $32$, $64$, $128$ and $256$
} 
\label{fig:Chentsov}
\end{figure}

\begin{figure}[htp]
    \begin{center}
    \includegraphics[width=5.5cm]{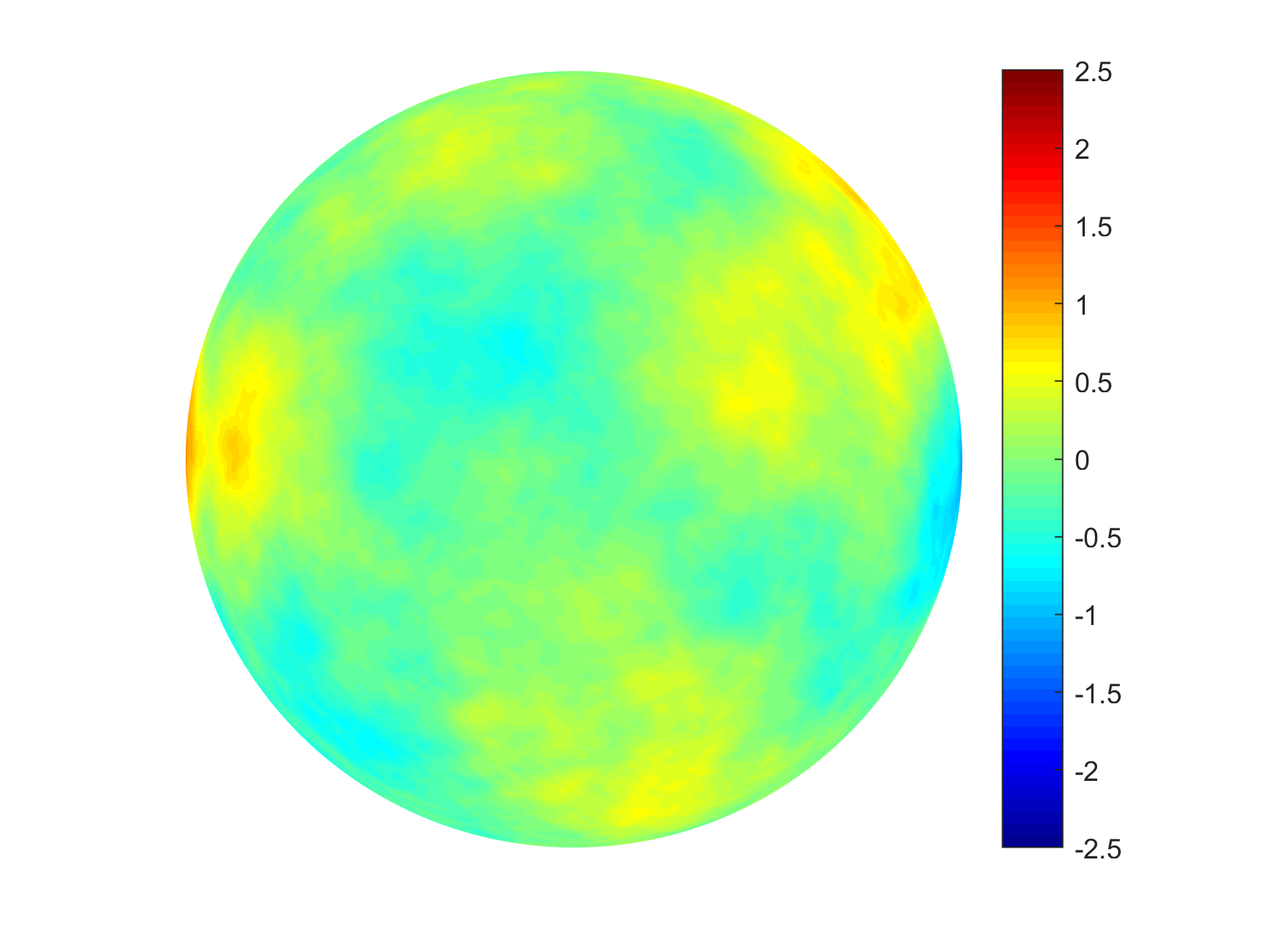} \end{center}
    \caption{Orthographic projection showing a realization of a univariate random field with a Chentsov covariance in $\mathbb{S}^{256}$, obtained by using $L=20,000$ basic random fields, all the other parameters being the same as in Figure \ref{fig:Chentsov}
} 
\label{fig:Chentsov20000}
\end{figure}

The isotropic Chentsov covariance on $\mathbb{S}^d$ is defined as $K_C(\vartheta) = 1-\frac{2\vartheta}{\pi}$ and its Schoenberg sequence $\{b_{n,d}^{C}: n\in\mathbb{N}\}$ is given in Appendix \ref{sub:sbg}. The following illustration (Figure \ref{fig:Chentsov}) displays orthographic projections of realizations on the $2$-sphere such that $x_1^2 + x_2^2 + x_3^2 = 1$ and $x_4 = \cdots = x_{d+1} = 0$ (intersection of $\mathbb{S}^d$ with the subspace whose last $d-2$ coordinates are zero), obtained with $L=1500$ basic random fields, a Rademacher distribution for $\varepsilon$, a zero probability for even integers $\kappa$ and a zeta distribution with parameter $2$ for odd integers $\kappa$, for dimensions $d$ ranging between $2$ and $256$. One notes that the striation effect is all the more pronounced as $d$ increases, which may be explained because the central limit approximation has a slower and slower rate of convergence. In particular, since the Schoenberg coefficient $b_{n,d}^{C}$ behaves like $n^{-d}$ as $n$ increases, the Berry-Ess\'{e}en bound is finite in the cases $d=2$ and $d=3$, but not necessarily for higher dimensions (Eq \eqref{thetaprimemax}). Interestingly, the striation effect becomes imperceptible when increasing the number of basic random fields to $L=20,000$ (Figure \ref{fig:Chentsov20000}).

\section{Practical aspects}
\label{practicalaspects}

\subsection{Distribution of $\kappa$}

The distribution of $\kappa$ should give a non-negligible probability to any degree having a significant contribution to the spectral representation of the target random field (degree $n$ for which the Schoenberg matrix $\bm{B}_{n,d}$ has large entries). In practice, many of the usual covariance models (with the exception of the negative binomial model) have a Schoenberg sequence that is lower bounded by a hyperharmonic series (behaving like $n^{-\theta}$ with $\theta>d-1$) and their rate of decay as $n$ increases is quite slow. Based on the third case presented in Section \ref{distributionchoice}, it is convenient to choose a zeta distribution for the random integer $\kappa$ (Eq. \eqref{thetaprimemax}) in order to ensure a finite Berry-Ess\'{e}en bound and a convergence to normality in $L^{-1/2}$. Such a distribution is long-tailed and allows the simulated random field to be a mixture of Gegenbauer waves with degrees ranging from very low to very high. This option, which has been adopted in Examples 2 to 4 above, is particularly interesting in order to reproduce both the low-frequency (large-scale) and high-frequency (small-scale) variations of the target random field.

However, when simulating on high-dimensional spheres or when the covariance model is associated with a Schoenberg sequence that is not lower bounded by a hyperharmonic series (which corresponds to a strongly irregular random field), the use of a zeta distribution for $\kappa$ may not guarantee the existence of a finite Berry-Ess\'{e}en bound. In such cases, one may trade the zeta distribution for a `super-heavy' tailed distribution, e.g., a distribution with a logarithmically decaying tail such as the discretized log-Cauchy distribution. The same issue arises with simulation algorithms based on spherical or hyperspherical harmonics approximations \citep{emery2019, lantuejoul2019}, with the inconvenience that the calculation of such harmonics for high degrees is particularly expensive and can make these algorithms prohibitive in terms of computation time. Also note that having a infinite Berry-Ess\'{e}en bound does not prevent the simulated random field to converge to a Gaussian random field as $L$ tends to infinity: it just means that the convergence rate can be slower than $L^{-1/2}$.

\subsection{Number of basic random fields (Gegenbauer waves)}

The choice of the number $L$ of basic random fields depends on the smoothness of the target random field and the dimension of the sphere on which it is simulated: as illustrated with the examples, more basic random fields are needed for irregular random fields (covariance function that quickly decays near the origin) and/or for high-dimensional spheres, in order to avoid the striation effect. The latter effect indicates that the convergence to multivariate normality is not reached, although the simulated field possesses the correct first- and second-order moments (expectation and covariance function). As a rule of thumbs, unless the target random field is strongly irregular or the simulation is performed on a high-dimensional sphere, a few thousand basic random fields ($L=1000$ to $5000$) is often sufficient to get `good-looking' realizations.

\subsection{Computer implementation and running time}

A set of Matlab\textsuperscript{\tiny\textregistered} scripts implementing turning arcs simulation is provided in Additional Material. These scripts consist of

\begin{itemize}
\item one main routine (\textit{turningarcs.m}) allowing the simulation of random fields on $\mathbb{S}^d$ with negative binomial, spectral-Matérn, generalized $\cal{F}$, Chentsov and exponential covariances, using a Rademacher distribution for $\varepsilon$ and a zeta distribution with parameter 2 for $\kappa$;

\item two subroutines (\textit{Gegenbauer.m} and \textit{zetarnd.m}) to calculate Gegenbauer polynomials based on \eqref{Gegenbauer} and to simulate $\kappa$, respectively;

\item one instruction file (\textit{examples.m}) that reproduces the examples shown in the Section \ref{examples}.

\end{itemize}

Executing the examples on a desktop with $128$ GB RAM and an Intel\textsuperscript{\tiny\textregistered} Xeon\textsuperscript{\tiny\textregistered} processor @$2.10$ GHz for simulating a random field on $\mathbb{S}^2$ discretized into $500 \times 500$ faces takes around $2$ seconds when using $L=15$ (Example 1 for a bivariate negative binomial covariance) and around $30$ seconds when using $L=1500$ (Examples 1, 2 and 4 for the negative binomial, spectral-Mat\'ern and Chentsov covariances). These running times, which include pre-processing, simulation and writing the results into an output ASCII file, are  smaller than that of algorithms based on spherical harmonics approximations, the most competitive ones being the spectral algorithms proposed by \citet{emery2019} and \citet{lantuejoul2019}. The difference in running time between the turning arcs and other spectral algorithms considerably increases when simulating at irregularly spaced colatitudes and longitudes, in which case the algorithms using spherical harmonics turn out to be much slower.

Still with $L = 1500$, the turning arcs algorithm takes $150$ seconds to simulate a random field with a generalized $\cal{F}$-covariance on $\mathbb{S}^3$ discretized into $8 \times 500 \times 500$ faces, from which the maps in Figure \ref{fig:Hypersphere} can be obtained (Example 3): the higher computation time ($5$ times more than for the examples in $\mathbb{S}^2$) is mainly explained because there are $8$ times more locations targeted for simulation in this example. As for Example 4 concerning the simulation of a random field with a Chentsov covariance on a sphere of dimension $256$ discretized into $500 \times 500 \times 1 \times \cdots \times 1$ faces, the computation time increases to $258$ seconds ($4.3$ minutes) with $L=1500$ basic random fields and $1110$ seconds ($18.5$ minutes) with $L=20,000$. All these examples prove that the simulation on spheres of more than $2$ dimensions is considerably faster than that based on expansions into hyperspherical harmonics, the computation of which is much more expensive than that of Gegenbauer polynomials.

Finally note that the turning arcs algorithm lends itself to parallel computing (not implemented in the Additional Material scripts), which could decrease all the aforementioned calculation times by one or two orders of magnitude.

\section{Conclusions}
\label{conclusions}

The turning arcs algorithm allows simulating isotropic scalar- and vector-valued Gaussian random fields on the sphere $\mathbb{S}^d$, provided that the spectral representation (Schoenberg sequence) of their covariance function is known. The simulation is obtained by spreading Gegenbauer waves that vary along randomly oriented arcs along the parallels orthogonal to these arcs, alike the continuous spectral and turning bands algorithms used to simulate random fields in Euclidean spaces. The advantages of the algorithm over existing alternatives are threefold: (1) it is extremely flexible, as it allows the simulation of vector random fields with any number of components, any isotropic covariance structure, on any $d$-dimensional sphere and any number and configuration of points targeted for simulation; (2) it accurately reproduces the desired covariance, and (3) it is computationally inexpensive, the numerical complexity being essentially proportional to the number of target points. Furthermore, with a suitable choice of the simulation parameters, the rate of convergence of the simulated random field to normality is at most of the order of $L^{-1/2}$, where $L$ is the number of Gegenbauer waves, except for covariance models on high-dimensional spheres ($d \geq 4$) whose Schoenberg coefficients decrease slowly as $n$ increases. A by-product of this research is a closed-form expression of the Schoenberg coefficients associated with the Chentsov and exponential covariance models in $\mathbb{S}^d$ for any $d \geq 2$.

\section*{Acknowledgements}

The authors acknowledge the funding of the National Agency for Research and Development of Chile, through grants CONICYT/FONDECYT/INICIACI\'ON/No.\ 11190686 (A. Alegr\'ia), CONICYT/FONDECYT/REGULAR/No.\ 1170290 (X. Emery) and CONICYT PIA AFB180004 (X. Emery).

%The research work of Alfredo Alegr\'ia was partially supported by grant CONICYT/FONDECYT/REGULAR/No.\ 11190686 from the Chilean government. Xavier Emery acknowledges the support of grants CONICYT/FONDECYT/REGULAR/No.\ 1170290 and CONICYT PIA AFB180004 (AMTC) from the Chilean Commission for Scientific and Technological Research. 

\appendix
\section*{Appendices}
\addcontentsline{toc}{section}{Appendices}
\renewcommand{\thesubsection}{\Alph{subsection}}
\renewcommand{\theequation}{\Alph{subsection}.\arabic{equation}}
\renewcommand\thetable{\thesubsection.\arabic{table}}

\subsection{Proof of Proposition \ref{prop1}}
\label{proof-prop1}

Before stating the proof of Proposition \ref{prop1}, we must introduce some properties of Gegenbauer polynomials. A classical duplication equation (see, e.g., \citealp[Equation 2.4]{ziegel2014convolution}) establishes that, for $d\geq 2$ and for any $n,k\in\mathbb{N}$,
\begin{equation}
\label{propiedad1}
\int_{\mathbb{S}^d}   {G}_n^{(d-1)/2}(\bm{\omega}^\top \bm{x}_1) \, {G}_k^{(d-1)/2}(\bm{\omega}^\top \bm{x}_2) \,  U({\rm d}\bm{\omega})   =  \frac{\delta_{n,k} (d-1)}{2n + d - 1} \, {G}_n^{(d-1)/2}(\bm{x}_1^\top \bm{x}_2),  \qquad \bm{x}_1,\bm{x}_2\in\mathbb{S}^d,
\end{equation}
where $U$ is the uniform probability measure on $\mathbb{S}^d$ and $\delta_{n,k}$ denotes the Kronecker delta.  %In particular, for any $n\in\mathbb{N}$, one has
%\begin{equation}
%\label{propiedad2}
%\int_{\mathbb{S}^d}    {G}_n^{(d-1)/2}(\bm{\omega}^\top \bm{x})  U({\rm d}\bm{\omega})   = 0,  \qquad \bm{x}\in\mathbb{S}^d.
%\end{equation}
For $d=1$, one has a similar identity. Let $n, k \in\mathbb{N}$, then
\begin{equation*}
\label{propiedad1_d=1}
\int_{\mathbb{S}^1}   \cos(n \vartheta(\bm{\omega}, \bm{x}_1)) \cos(k \vartheta(\bm{\omega}, \bm{x}_2)) \, U({\rm d}\bm{\omega})   =  \frac{\delta_{n,k}}{2} \cos(n \vartheta(\bm{x}_1, \bm{x}_2)),  \qquad \bm{x}_1,\bm{x}_2\in\mathbb{S}^1.
\end{equation*}
%Again, as a particular case, one has that for any $n\in\mathbb{N}$,
%\begin{equation*}
%\label{propiedad2_d=1}
%\int_{\mathbb{S}^1}   \cos(n \vartheta(\bm{\omega}, \bm{x}))  U({\rm %d}\bm{\omega})   =  0,  \qquad \bm{x}\in\mathbb{S}^1.
%\end{equation*}
%In addition, $$\int_{\mathbb{S}^1}   {G}_0^{0}(\bm{\omega}^\top \bm{x}_1) {G}_0^{0}(\bm{\omega}^\top \bm{x}_2)  U({\rm d}\bm{\omega})   =  1,\qquad  \bm{x}_1,\bm{x}_2\in\mathbb{S}^1.$$

\textbf{Proof of Proposition \ref{prop1}}
We only prove the result for $d\geq 2$, since the case $d=1$ is completely analogous. Let ${Z}$ be the random field defined in (\ref{simulacion2}). Because $\varepsilon$ is independent of $\kappa$ and $\bm{\omega}$ and has a zero mean, it is straightforward to prove that $\mathbb{E}\{{Z}(\bm{x})\} = 0$ for any $\bm{x}\in\mathbb{S}^d$.
%The mean is given by
%\begin{equation}
%\label{media2}
%\mathbb{E}\{{Z}(\bm{x})\}    =     \sum_{n=0}^\infty   a_n  \sqrt{\frac{b_{n,d} (2n + d - 1)}{a_n (d-1)} }  \int_{\mathbb{S}^d}  {G}_n^{(d-1)/2}(\bm{\omega}^\top \bm{x})   U(\text{d}{\bm{\omega}}), \qquad \bm{x}\in\mathbb{S}^d.
% \end{equation}
%Thus, the mean is obtained from property (\ref{propiedad2}), since it implies that all terms in (\ref{media2}) are identically equal to $0$, except for $n=0$, for which the last integral is identically equal to one.  
On the other hand, the covariance between any two variables ${Z}(\bm{x}_1)$ and ${Z}(\bm{x}_2)$, with $\bm{x}_1,\bm{x}_2\in\mathbb{S}^d$, is:
\begin{equation*}
\mathbb{E}\{{Z}(\bm{x}_1) {Z}(\bm{x}_2)\}     =  \mathbb{E}\{\varepsilon^2\}  \sum_{n=0}^\infty \frac{b_{n,d} (2n+d-1)}{d-1} \int_{\mathbb{S}^d}  {G}_n^{(d-1)/2}(\bm{\omega}^\top \bm{x}_1) \, {G}_n^{(d-1)/2}(\bm{\omega}^\top \bm{x}_2) \, U(\text{d}{\bm{\omega}}). %, \qquad \bm{x}_1,\bm{x}_2\in\mathbb{S}^d.
 \end{equation*}
Using (\ref{propiedad1}) and the fact that $\mathbb{E}\{\varepsilon^2\}=1$, the announced covariance function is obtained.

\subsection{Proof of Proposition \ref{prop2}}
\label{proof-prop2}

Again, we only prove the result for $d\geq 2$, the one-dimensional case being similar. Let $\bm{Z}$ be the vector-valued random field defined in (\ref{simulacion-vector}). Its mean vector is zero, insofar as $\varepsilon$ has zero mean and is independent of $\bm{\omega}$, $\iota$ and $\kappa$. 
%given by
%\begin{equation}
%\label{media2-vector}
%\mathbb{E}\{\bm{Z}(\bm{x})\}    =    \sum_{i=1}^p \sum_{n=0}^\infty a_n  \sqrt{ \frac{ 2n + d -1}{a_{n} (d-1) } } \bm{\gamma}^{(i)}_{n,d}    \int_{\mathbb{S}^d}  {G}_n^{(d-1)/2}(\bm{\omega}^\top \bm{x})   U(\text{d}{\bm{\omega}}), \qquad \bm{x}\in\mathbb{S}^d.
% \end{equation}
%One obtains the mean vector by using property (\ref{propiedad2}). Indeed, it implies that all terms in (\ref{media2-vector}) vanish, except for $n=0$.  

The variance-covariance matrix between any two vectors $\bm{Z}(\bm{x}_1)$ and $\bm{Z}(\bm{x}_2)$, with $\bm{x}_1,\bm{x}_2\in\mathbb{S}^d$, is:
%In order to obtain the covariance function, note that
\begin{equation*}
\begin{split}
%\begin{multline*}
\mathbb{E} & \{\bm{Z}(\bm{x}_1) \bm{Z}(\bm{x}_2)^\top\}  \nonumber \\
  & =  \mathbb{E} \{\varepsilon^2 \} \sum_{n=0}^\infty a_n  \left\{ \frac{p (2n + d -1)}{a_n (d-1) } \right\}  \frac{1}{p} \sum_{i=1}^p \bm{\gamma}^{(i)}_{n,d}\, [\bm{\gamma}^{(i)}_{n,d}]^\top \int_{\mathbb{S}^d}   {G}_{n}^{(d-1)/2}(\bm{\omega}^\top \bm{x}_1) \,  {G}_{n}^{(d-1)/2}(\bm{\omega}^\top \bm{x}_2)   \,  U(\text{d}{\bm{\omega}}). \nonumber  \\
% \end{multline*}
\end{split}
\end{equation*}
 Using property (\ref{propiedad1}) and the fact that $\varepsilon$ is an independent random variable with zero mean and unit variance, one obtains
 %(\ref{propiedad2}), one obtains 
$$ \mathbb{E}\{\bm{Z}(\bm{x}_1) \bm{Z}(\bm{x}_2)^\top\}     =  \sum_{n=0}^\infty   \left\{\sum_{i=1}^p \bm{\gamma}^{(i)}_{n,d}[\bm{\gamma}^{(i)}_{n,d}]^\top\right\}  G_n^{(d-1)/2}(\bm{x}_1^\top \bm{x}_2). $$  %+  a_0 \sum_{\substack{i,j=1 \\ j\neq i}}^p  \bm{\gamma}_{0,d}^{(i)}[\bm{\gamma}_{0,d}^{(j)}]^\top.$$
%Also, note that 
%$$  \mathbb{E}\{\bm{Z}(\bm{x}_1)\} \mathbb{E}\{ \bm{Z}(\bm{x}_2)\}^\top  = a_0 \sum_{i,j=1}^p   \bm{\gamma}^{(i)}_{0,d}[\bm{\gamma}^{(j)}_{0,d}]^\top.     $$ 
%One concludes that 
%$$ \mathbb{E}\{\bm{Z}(\bm{x}_1) \bm{Z}(\bm{x}_2)^\top\}  -  \mathbb{E}\{\bm{Z}(\bm{x}_1)\} \mathbb{E}\{ \bm{Z}(\bm{x}_2)\}^\top   =  \sum_{n=0}^\infty   \left\{\sum_{i=1}^p \bm{\gamma}^{(i)}_{n,d}[\bm{\gamma}^{(i)}_{n,d}]^\top\right\}  G_n^{(d-1)/2}(\bm{x}_1^\top \bm{x}_2)  - a_0 \sum_{i=1}^p   \bm{\gamma}^{(i)}_{0,d}[\bm{\gamma}^{(i)}_{0,d}]^\top.$$
The covariance function is obtained by using (\ref{gamma-matrix}).

\subsection{Upper bound for the third-order absolute moment of a Gegenbauer wave}
\label{proof-BerryEsseen}

Let $d, n \in \mathbb{N}$, $d \geq 2$, $\lambda = \frac{d-1}{2}$, $\bm{x} \in \mathbb{S}^d$ (fixed) and $\bm{\omega}$ uniformly distributed on $\mathbb{S}^d$. It is of interest to find an upper bound for the following third-order absolute moment:
$$ \mu_{n,d}^3 = \mathbb{E} \{ \lvert G_n^\lambda(\bm{\omega}^T \, \bm{x}) \rvert^3 \}. $$

By introducing spherical coordinates such that:
$$
\left\{
\begin{aligned} 
       \bm{x} &=  (1,0,\cdots,0) \\
      \bm{\omega} &=  (\cos \varphi_1,\sin \varphi_1 \cos \varphi_2, \cdots,\sin \varphi_1 \cdots \sin \varphi_{d-1} \cos \varphi_d,\sin \varphi_1 \cdots \sin \varphi_{d-1} \sin \varphi_d),
\end{aligned}
\right.
$$
with $\varphi_1, \cdots, \varphi_{d-1} \in [0,\pi]$ and $\varphi_d \in [0,2\pi[$, one obtains:
$$ 
\begin{aligned} 
\mu_{n,d}^3 &= \int_{\mathbb{S}^d} \lvert G_n^\lambda(\bm{\omega}^T \, \bm{x}) \rvert^3 U({\rm d}\bm{\omega}) \\
 &= \frac{\Gamma\left(\frac{d+1}{2}\right)}{2\pi^{\frac{d+1}{2}}} \int_0^{2\pi} {\rm d}\varphi_d \int_0^{\pi} \sin \varphi_{d-1} {\rm d}\varphi_{d-1} \cdots \int_0^{\pi} \sin^{d-2} \varphi_{2} {\rm d}\varphi_{2} \int_0^{\pi} \lvert G_n^\lambda(\cos \varphi_{1}) \rvert^3 \sin^{d-1} \varphi_{1} {\rm d}\varphi_{1}.
\end{aligned}
$$
\\
Since $\int_0^{\pi} \sin^{m-1} \varphi = \frac{\sqrt{\pi} \Gamma(\frac{m}{2})}{\Gamma(\frac{m+1}{2})}$ (\citet{grad}, formula 3.621.5), one has: 
\begin{equation}
\label{mund3}
    \mu_{n,d}^3 = \frac{2\Gamma\left(\frac{d+1}{2}\right)}{\sqrt{\pi}\Gamma\left(\frac{d}{2}\right)} \int_0^{\frac{\pi}{2}} \lvert G_n^\lambda(\cos \varphi_{1}) \rvert^3 \sin^{d-1} \varphi_{1} {\rm d}\varphi_{1}.
\end{equation}

To find an upper bound for such a moment, we distinguish the cases $d=2$, $d=3$ and $d \geq 4$.
\begin{itemize}
    \item [Case $d=2$.] Using inequality 22.14.3 of \citep{abramowitz}:
    $$\lvert G_n^{\frac{1}{2}}(\cos \varphi_{1}) \rvert \leq \sqrt{\frac{2}{n \pi \sin \varphi_1}}, \qquad \varphi_1 \in ]0,\pi[,$$
    one finds
\begin{equation*}
    \mu_{n,2}^3 \leq \frac{2\Gamma\left(\frac{d+1}{2}\right)}{\sqrt{\pi}\Gamma\left( \frac{d}{2}\right)} \left( \frac{2}{n \pi}\right)^{\frac{3}{2}} \int_0^{\frac{\pi}{2}} \frac{1}{\sqrt{\sin(\varphi_{1})}} {\rm d}\varphi_{1},
\end{equation*}
i.e.,
\begin{equation}
    \mu_{n,2}^3 \leq \left(\frac{2}{n \pi}\right)^{\frac{3}{2}} \frac{\Gamma\left(\frac{1}{4}\right)}{\pi\Gamma\left( \frac{3}{4}\right)} = \mathcal{O}\left(n^{-\frac{3}{2}}\right).
\end{equation}

    \item [Case $d=3$.] In this case, the Gegenbauer polynomials of order $\lambda = 1$ coincide with the Chebyshev polynomials of the second kind (\citet{abramowitz}, formulae 22.5.34 and 22.3.16): 
    $$
    G_n^1(\cos \varphi_1) = \frac{\sin\left((n+1)\varphi_1 \right)}{\sin (\varphi_1)}, \qquad n \in \mathbb{N},  \varphi_1 \in ]0,\pi[.
    $$
    We use the following inequalities:
$$ 
\left\{
\begin{aligned} 
\lvert \sin\left((n+1)\varphi_1\right) \rvert &\leq (n+1)\varphi_1 \text{ for }  0\leq \varphi_1 \leq \frac{\pi}{2(n+1)} \\
\lvert \sin\left((n+1)\varphi_1\right) \rvert &\leq 1 \text{ for } \frac{\pi}{2(n+1)} \leq \varphi_1 \leq \frac{\pi}{2} \\
\lvert \sin(\varphi_1) \rvert &\geq \frac{2}{\pi}\varphi_1 \text{ for } 0\leq \varphi_1 \leq \frac{\pi}{2},
\end{aligned} 
\right.
$$
which yield:
\begin{equation*}
    \mu_{n,3}^3 \leq \frac{2\Gamma\left(\frac{d+1}{2}\right)}{\sqrt{\pi}\Gamma\left(\frac{d}{2}\right)} \left( \frac{\pi^3 (n+1)^3}{2^3}  \int_0^{\frac{\pi}{2(n+1)}} \varphi_{1}^2 {\rm d}\varphi_{1} + \frac{\pi^3}{2^3}  \int_{\frac{\pi}{2(n+1)}}^{\frac{\pi}{2}}  \frac{{\rm d}\varphi_{1}}{\varphi_{1}} \right),
\end{equation*}
that is:
\begin{equation}
    \mu_{n,3}^3 \leq \frac{\pi^5}{48} + \frac{\pi^2 \ln(n+1)}{2} = \mathcal{O}\left(\ln n\right).
\end{equation}

    \item [Case $d \geq 4$.] Let us pose $\lambda = \frac{d-1}{2}$. For any integer $\nu \in [1,\lambda[$, \citet{Reimer} showed that there exists a constant $\varpi_{\nu,n}^{\lambda}$ depending on $\nu$, $n$ and $\lambda$ such that
    %One has \citep{Reimer}:
    $$
    \lvert G_n^{\lambda}(\cos \varphi_{1}) \rvert \leq \varpi_{\nu,n}^{\lambda} G_n^{\lambda}(1) \lvert n \sin(\varphi_1) \rvert ^{-\nu}, \qquad \varphi_1 \in ]0,\pi[.
    $$
%where $\nu$ can be any integer such that $\lambda > \nu \geq 1$ and $\varpi_{\nu,n}^{\lambda}$ is a constant depending on $\nu$, $n$ and $\lambda = \frac{d-1}{2}$. 
Plugging this inequality into \eqref{mund3}, one obtains:
    $$
    \mu_{n,d}^3 \leq \frac{2\Gamma\left(\frac{d+1}{2}\right)}{\sqrt{\pi}\Gamma\left(\frac{d}{2}\right)} \left(\varpi_{\nu,n}^{\lambda} \, G_n^{\lambda}(1)\right)^3 \, n^{-3\nu} \int_0^{\frac{\pi}{2}} \sin^{d-1-3\nu}(\varphi_1) {\rm d}\varphi_{1},
    $$
with $G_n^{\lambda}(1) = \frac{\Gamma(n+2\lambda)}{\Gamma(2\lambda)\Gamma(n+1)}$ (\citet{abramowitz}, formula 22.2.3). 

The above integral converges when $d-1-3\nu$ is greater than $-1$ (\citet{grad}, formula 3.621.5), in which case one has:
    $$
    \mu_{n,d}^3 \leq \frac{2\Gamma\left(\frac{d+1}{2}\right)}{\sqrt{\pi}\Gamma\left(\frac{d}{2}\right)} \left(\varpi_{\nu,n}^{\lambda} \, \frac{\Gamma(n+d-1)}{\Gamma(d-1)\Gamma(n+1)}\right)^3 \, n^{-3\nu} \frac{\sqrt{\pi}\Gamma \left( \frac{d-3\nu}{2} \right)}{2\Gamma \left( \frac{d-3\nu+1}{2} \right)}.
    $$

\cite{Reimer} showed that $\varpi_{\nu,n}^{\lambda} = \mathcal{O}(1)$ as $n$ becomes infinitely large. Furthermore, Stirling's approximation to the factorial implies that $\frac{\Gamma(n+d-1)}{\Gamma(n+1)} = \mathcal{O}\left(n^{d-2}\right)$ (\citet{abramowitz}, formula 6.1.46). The lowest asymptotic bound is obtained for $\nu=\lfloor \frac{d}{2} \rfloor - 1$:
\begin{equation}
    \mu_{n,d}^3 \leq \mathcal{O}\left(n^{3d-6-3\nu}\right) = \mathcal{O}\left(n^{3 \lfloor \frac{d-1}{2} \rfloor}\right).
\end{equation}

\end{itemize}

%%%%%%%%%%%%%%%%%%%%%%%%%%%%%%%%%%%%%%%%%%%%%%%%%%%%%%%%%%%%%%%%%%%%%%%%%%%%%%%%%%%%%%

%
\subsection{Calculation of Schoenberg coefficients}\label{sub:sbg} 
Let $K$ be an isotropic, positive semi-definite function on $\mathbb{S}^d$, $d \geq 2$, and let $\{b_{n,d}: n\in\mathbb{N}\}$ be its associated Schoenberg sequence, as defined in \eqref{schoenberg1}. The change of variable $t = \cos \vartheta$ in \eqref{bnd1} gives

%recall its series representation \eqref{schoenberg1} into Gegenbauer polynomials: 
%
%$$ K (\vartheta) = \sum_{n=0}^\infty b_{n,d} \, G_n^\lambda ( \cos \vartheta) $$
%
%with $\lambda = \frac{d - 1}{2}$. The Schoenberg coefficients $\{b_{n,d}: n\in\mathbb{N}\}$ can be retrieved using the formula 
%
%$$ b_{n,d} = \frac{1}{\parallel G_n^\lambda \parallel^2} \int_0^\pi K (\vartheta) \, G_n^\lambda ( \cos \vartheta) \, \sin^{ 2 \lambda} \, d \vartheta $$
%
%with 
%
%$$ \parallel G_n^\lambda \parallel^2 = \int_0^\pi \bigl[ G_n^\lambda (\cos \vartheta) \bigr]^2 \, \sin^{2 \lambda} (\vartheta) \, d \vartheta. $$
%
%The change of variable $t = \cos \vartheta$ in \eqref{bnd1} gives
%
$$ b_{n,d} = \frac{1}{\parallel G_n^\lambda \parallel^2} \int_{-1}^{+1} K (\arccos t) \, G_n^\lambda (t) \, 
\bigl(1 - t^2 \bigl)^{\lambda - 1/2} \, d t, $$

with $\lambda = \frac{d - 1}{2} > 0$. Suppose now that $t \, \mapsto \, K (\arccos t)$ can be expanded into a power series
$$ K (\arccos t) = \sum_{k=0}^\infty \alpha_k \, t^k. $$
Then, using the expansion of the monomials into Gegenbauer polynomials \citep{rve60, Kim} 
$$ t^k = \frac{k!}{2^k} \, \sum_{\ell=0}^{\lfloor k/2 \rfloor} \frac{\lambda + k - 2 \ell}{\ell !} \, \frac{\Gamma (\lambda)}
{\Gamma ( \lambda + k + 1 -\ell)} \, G_{k - 2 \ell}^\lambda (t)  , $$
where $\lfloor \cdot \rfloor$ is the floor function, it follows
$$ b_{n,d} = \frac{1}{\parallel G_n^\lambda \parallel^2} \sum_{k=0}^\infty \alpha_k \, \frac{k!}{2^k}  \sum_{\ell=0}^{\lfloor k/2 \rfloor} 
\frac{\lambda + k - 2 \ell}{\ell !} \, \frac{\Gamma (\lambda)}{\Gamma ( \lambda + k + 1 -\ell)} \,
\int_{-1}^{+1} G_{k - 2 \ell}^\lambda (t) \, G_n^\lambda (t) \, \bigl(1 - t^2 \bigl)^{\lambda - 1/2} \, d t. $$
The latter integral vanishes unless $k - 2 \ell = n $, in which case it is equal to $\parallel G_n^\lambda \parallel^2$.
%
%$$ b_{n,d} = \frac{1}{\parallel G_n^\lambda \parallel^2} \sum_{k=0}^\infty \alpha_k \sum_{\ell=0}^{\lfloor k/2 \rfloor} \frac{\lambda + k - 2 \ell}{\ell !} \, \frac{\Gamma (\lambda)}{\Gamma ( \lambda + k + 1 -\ell)} \, \parallel G_n^\lambda \parallel^2 \, 1_{k - 2 \ell = n}. $$
%
We thus obtain the generic formula 
\begin{equation}\label{eq:bnd}
b_{n,d} = \sum_{\ell=0}^\infty \alpha_{n + 2 \ell} \, \frac{(n + 2 \ell)!}{2^{n + 2 \ell}} \, \frac{\lambda + n}{\ell !} 
\, \frac{\Gamma (\lambda)}{\Gamma ( \lambda + n + \ell + 1)} \cdot
\end{equation}
The rest of the calculation must be done on a case-by-case basis. Two examples are given below. 

\subsubsection{Chentsov covariance}
As a first example, consider $K (\vartheta) = 1 - 2 \vartheta / \pi$. The power series of $ K (\arccos t) = \frac{2}{\pi} \, \arcsin t $ is given  
by formula 4.4.40 of \citet{abramowitz}:  
$$ K ( \arccos t) = \frac{2}{\pi \sqrt{\pi}} \, \sum_{k=0}^\infty \frac{\Gamma (k + 1/2)}{(2 k +1) \, k!} \, t^{2k+1} , $$
from which we derive $\alpha_{2k} = 0$ and
$$ \alpha_{2k+1} = \frac{2}{\pi \sqrt{\pi}} \, \frac{\Gamma (k + 1/2)}{(2 k +1) \, k!} \cdot $$ 
Plugging these coefficients into \eqref{eq:bnd}, we obtain that $b_{2n,d}= 0$ and
$$ b_{2n+1,d} = \frac{2}{\pi \sqrt{\pi}} \sum_{\ell=0}^\infty \frac{\Gamma ( n + \ell + 1/2)}{(2n+2\ell+1) \, 
( n + \ell)!} \, \, \frac{(2n + 2 \ell +1)!}{2^{2n + 2 \ell+1}} \, \frac{\lambda + 2n+1}{\ell !} \, \frac{\Gamma (\lambda)}
{\Gamma ( \lambda + 2n + \ell + 2)}. $$
Using the duplication formula of the gamma function (formula 6.1.18 of \citet{abramowitz}), it comes
\begin{equation*}
\begin{split}
b_{2n+1,d} &= \frac{1}{\pi^2} \sum_{\ell=0}^\infty \frac{\Gamma^2 ( n + \ell + 1/2) \, (\lambda + 2n+1) \, \Gamma (\lambda)}{\ell ! \, \Gamma ( \lambda + 2n + \ell + 2)} \\
&= \frac{(\lambda + 2n+1) \, \Gamma (\lambda)}{\pi^2} \, \frac{\Gamma^2 ( n + 1/2)}{\Gamma ( \lambda +2n+2)} \, {}_2F_1 \left(n+\frac{1}{2}, n+\frac{1}{2}; \lambda +2n+2; 1\right),
\end{split}
\end{equation*}
where ${}_2F_1$ is the Gaussian hypergeometric function.
%The summation can be expressed with the Gaussian hypergeometric function ${}_2F_1$
%leads to an hypergeometric function
%
%$$ b_{2n+1,d} = \frac{(\lambda + 2n+1) \, \Gamma (\lambda)}{\pi^2} \, \frac{\Gamma^2 ( n + 1/2)}{\Gamma ( \lambda +2n+2)} \, {}_2F_1 (n+1/2, n+1/2; \lambda +2n+2; 1) , $$
Owing to Gauss's theorem (formula 15.1.20 of \citet{abramowitz}), this finally reduces to 
\begin{equation}\label{eq:csvbnd}
b_{2n+1,d} = \frac{(\lambda + 2n+1) \, \Gamma (\lambda) \, \Gamma (\lambda +1)}{\pi^2} \, \frac{\Gamma^2 ( n + 1/2)}
{\Gamma^2 ( \lambda + n + 3/2)}.    
\end{equation}
%
%by applying Gauss' theorem (formula 15.1.20 of \citet{abramowitz}).
These coefficients can be calculated directly, or by using the 
induction formula
$$ b_{2 n+1,d} = \frac{\lambda + 2 n +1}{\lambda + 2 n -1} \, \frac{(n-1/2)^2}{(\lambda + n + 1/2)^2} \, 
b_{2 n -1,d}, \qquad n \geq 1, $$
starting from
$$ b_{1,d} = \frac{\Gamma (\lambda) \, \Gamma (\lambda +2)}{\pi \, \Gamma^2 (\lambda +3/2)}. $$
 Equation \eqref{eq:csvbnd} generalizes the expressions provided by \citet{hagzhgrbn11} and \citet{lantuejoul2019} for the specific case when $d=2$. 

\subsubsection{Exponential covariance}
Put $K ( \vartheta ) = \exp ( - \nu \vartheta )$ with $\nu>0$. The power series of $ F(t) = \exp ( - \nu \arccos t )$ is 
required. A first derivation gives $ \sqrt{1 - t^2} F^\prime (x) - \nu F (t) = 0$. A second derivation followed by a multiplication 
by $\sqrt{ 1 - t^2}$ yields $ (1-t^2) F^{\prime\prime} (t) - t F^\prime (t) - \nu \sqrt{1 - t^2} F^\prime (t) = 0$. 
Replacing the third term by its expression in the first derivation, we finally obtain 
\begin{equation}\label{eq:dif}
( 1 - t^2) F^{\prime\prime} (t) - t F^\prime (t) - \nu^2 F (t) = 0. 
\end{equation}
Let us now expand $F$ into a power series:
$$ F (t) = \sum_{k=0}^\infty \alpha_k t^k.$$
Owing to the expression of $F$ and to the first derivation formula, the first two coefficients are $\alpha_0 = \exp (- \nu \pi /2)$ 
and $ \alpha_1 = \nu \exp ( - \nu \pi / 2)$. More generally, if the the power series of $F^\prime$ and $F^{\prime\prime}$ 
are plugged into \eqref{eq:dif}, then the following induction formula is obtained:
$$ \alpha_{k+2} = \alpha_k \, \frac{k^2 + \nu^2}{(k+1)(k+2)} = \alpha_k \, \frac{4}{(k+1)(k+2)} \left( \frac{k+i\nu}{2}\right) \left(\frac{k-i\nu}{2}\right),$$
where $i$ is the imaginary unit. If $k$ is even, we have
%For even indices, we have
%
%\begin{align*}
%\alpha_{2 k}&= \frac{\alpha_0}{(2 k)!} \prod_{\ell=0}^{k-1} ( 4 \ell^2 + \nu^2) \\
%&= \alpha_0 \, \frac{2^{2 k}}{(2 k)!} \, \prod_{\ell=0}^{k-1} ( \ell + i \nu / 2 )  ( \ell - i \nu / 2 ) \\
$$\alpha_{k} = \alpha_0 \, \frac{2^{k}}{k!} \, \frac{\Gamma ( \frac{k + i \nu}{2} )}{\Gamma (\frac{i \nu}{2} )} \, \frac{\Gamma ( \frac{k - i \nu}{2} )}{\Gamma ( - \frac{i \nu}{ 2} )} = \alpha_0 \, \frac{2^{k}}{k!} \, \left\lvert \Gamma \left( \frac{k + i \nu}{2} \right) \right\rvert^2 \, 
\frac{\nu \sinh (\frac{\pi \nu}{2})}{ 2 \pi},$$ %\qquad \text{(formula 6.1.29 of \citet{abramowitz})} \\
%&= \frac{\nu \bigl( 1 - \exp ( - \pi \nu) \bigr)}{4 \pi} \frac{2^{k}}{k!} \, \left\lvert \Gamma \left( \frac{k + i \nu}{2} \right) \right\rvert^2, 
%\end{align*}
the last equality being obtained by using formula 6.1.29 of \citet{abramowitz}. Likewise, if $k$ is odd, we have
%For odd indices, we have
%
%\begin{align*}
%\alpha_{2 k + 1} &= \frac{\alpha_1}{(2k+1)!} \prod_{\ell=0}^{k-1} \bigl( (2 \ell + 1)^2 + \nu^2 \bigr) \\
%&= \alpha_1 \, \frac{2^{2 k}}{(2 k+1)!} \, \prod_{\ell=0}^{k-1} ( \ell + 1/2 + i \nu / 2 )  ( \ell + 1/2 - i \nu / 2 ) \\
$$\alpha_{k} = \alpha_1 \, \frac{2^{k-1}}{k!} \, \frac{\Gamma ( \frac{k + i \nu}{2})}{\Gamma (\frac{1 + i \nu}{2})} \, 
\frac{\Gamma ( \frac{k - i \nu}{2})}{\Gamma (\frac{1 - i \nu}{2})} = \alpha_1 \, \frac{2^{k}}{k!} \, \left\lvert \Gamma \left( \frac{k + i \nu}{2} \right) \right\rvert^2 \, 
\frac{\cosh (\frac{\pi \nu}{2})}{2 \pi},$$ \\%\qquad \text{(formula 6.1.30 of \citet{abramowitz})} \\
%&= \frac{\nu \bigl( 1 + \exp ( - \pi \nu) \bigr)}{4 \pi} \, \frac{2^{k}}{k!} \, \left\lvert \Gamma \left( \frac{k + i \nu}{2} \right) \right\rvert^2, 
%\end{align*}
%
the last equality being obtained by using formula 6.1.30 of \citet{abramowitz}.

\medskip

Accordingly, accounting for the above expressions of $\alpha_0$ and $\alpha_1$, in all cases we have
\begin{equation}\label{eq:expak}
\alpha_k = C (\nu, k) \, \frac{2^k}{k!} \, \left\lvert \Gamma \left( \frac{ k + i \nu}{2} \right) \right\rvert^2, 
\end{equation}
where
$$ C(\nu,k) = 
\begin{cases}
\displaystyle \frac{\nu \, \exp( - \pi \nu / 2) \, \sinh ( \pi \nu / 2)}{2 \pi} & \text{if $k$ is even} \\[0.2cm]
\displaystyle \frac{\nu \, \exp( - \pi \nu / 2) \, \cosh ( \pi \nu / 2)}{2 \pi} & \text{if $k$ is odd.}
\end{cases}
$$
Plugging this expression into formula \eqref{eq:bnd}, we obtain
%
%$$ b_{n,d} = \frac{\nu \bigl( 1 - (-1)^k \exp ( - \pi \nu) \bigr)}{4 \pi} \, \sum_{\ell=0}^\infty \Gamma \left( \ell + \frac{ n + i \nu}{2} \right) \,  \Gamma \left( \ell + \frac{ n - i \nu}{2} \right) \, \frac{\lambda + n}{\ell !} \, \frac{\Gamma (\lambda)}{\Gamma ( \lambda + n + \ell + 1)} $$
%
%
%$$ b_{n,d} = \frac{\nu \bigl( 1 - (-1)^k \exp ( - \pi \nu) \bigr)}{4 \pi} \, \sum_{\ell=0}^\infty \Gamma \bigl( \ell + (n + i \nu) /2 \bigr) \,  \Gamma \bigl( \ell + (n - i \nu) / 2 \bigr) \, \frac{\lambda + n}{\ell !} \, \frac{\Gamma (\lambda)}{\Gamma ( \lambda + n + \ell + 1)} $$
%
%
%$$ b_{n,d} = (\lambda + n) \, \Gamma (\lambda) \, \frac{\nu \bigl( 1 - (-1)^k \exp ( - \pi \nu) \bigr)}{4 \pi} \, \sum_{\ell=0}^\infty \frac{\Gamma \bigl( \ell + (n + i \nu) /2 \bigr) \,  \Gamma \bigl( \ell + (n - i \nu) / 2 \bigr)} {\ell ! \, \Gamma ( \lambda + n + \ell + 1)} $$
%
\begin{align*}
b_{n,d} &= C (\nu,n) \, \sum_{\ell=0}^\infty \Gamma \left( \ell + \frac{ n + i \nu}{2} \right) \, \Gamma \left( \ell + 
\frac{ n - i \nu}{2} \right) \, \frac{\lambda + n}{\ell !} \, \frac{\Gamma (\lambda)}{\Gamma ( \lambda + n + \ell + 1)} \\
& = (\lambda + n) \, \Gamma (\lambda) \, C ( \nu,n) \, \frac{\Gamma \bigl( \frac{n + i \nu}{2} \bigr) \, \Gamma 
\bigl( \frac{n - i \nu}{2} \bigr)}{ \Gamma ( \lambda + n + 1)} \, {}_2F_1 \left( \frac{n + i \nu}{2} , \frac{n - i \nu}{2} ; \lambda + n + 1 ; 1 \right).
\end{align*}
%where ${}_2F_1$ is the Gaussian hypergeometric function. By Gauss's theorem (\citet{abramowitz}, formula 15.1.20), it comes 
%This can be written in terms of the Gaussian hypergeometric function ${}_2F_1$: 
%
%$$ b_{n,d} = (\lambda + n) \, \Gamma (\lambda) \, \frac{\nu \bigl( 1 - (-1)^k \exp ( - \pi \nu) \bigr)}{4 \pi} \, \frac{\Gamma \bigl( \frac{n + i \nu}{2} \bigr) \,  \Gamma \bigl( \frac{n - i \nu}{2} \bigr)}{ \Gamma ( \lambda + n + 1)} \, {}_2F_1 \left( \frac{n + i \nu}{2} , \frac{n - i \nu}{2} ; \lambda + n + 1 ; 1 \right) $$
%
By Gauss's theorem, it comes 
%
%$$ b_{n,d} = (\lambda + n) \, \Gamma (\lambda) \, \frac{\nu \bigl( 1 - (-1)^k \exp ( - \pi \nu) \bigr)}{4 \pi} \, \frac{\Gamma \bigl( (n + i \nu) /2 \bigr) \,  \Gamma \bigl( (n - i \nu) / 2 \bigr)}{ \Gamma ( \lambda + n + 1)} \, \frac{ \Gamma ( \lambda + n + 1) \, \Gamma (\lambda +1)}{\Gamma \bigl( \lambda +1 + (n + i \nu) /2 \bigr) \, \lambda +1 + (n - i \nu) /2 \bigr)} $$
\begin{equation}\label{eq:expbnd}
b_{n,d} = C (\nu,n) \, (\lambda + n) \, \Gamma (\lambda) \, \Gamma (\lambda +1) \, \frac{ \displaystyle \left\lvert 
\Gamma \left(\frac{n + i \nu}{2} \right) \right\rvert^2}{ \displaystyle \left\lvert \Gamma \left( \lambda +1 + \frac{n + 
i \nu}{2} \right) \right\rvert^2}. 
\end{equation}
%
%\smallskip

%with $ \lambda = ( d - 1) /2$. 

Calculating the squared modulus of the complex-valued gamma function in the numerator of \eqref{eq:expbnd} can be done 
by applying the induction formula
%Regarding the calculation of the complex-valued gamma functions, it can be noted that the real part of their argument is integer or half integer. Accordingly, the induction formula  
%
$$ \left\lvert \Gamma \left( \frac{n + i \nu}{2} \right) \right\rvert^2 = \frac{ (n-2)^2 + \nu^2}{4} \, 
\left\lvert \Gamma \left( \frac{n-2  + i \nu}{2} \right) \right\rvert^2, \qquad n \geq 2, $$
along with the initial values (\citet{abramowitz}, formulae 6.1.29 and 6.1.30)
$$ \left\lvert \Gamma \left( \frac{i \nu}{2} \right) \right\rvert^2 = \frac{ 2 \pi}{\nu \, \sinh (\pi \nu / 2)} 
\qquad  \qquad \left\lvert \Gamma \left( \frac{1 + i \nu}{2} \right) \right\rvert^2 = \frac{\pi}{\cosh 
( \pi \nu / 2)}. $$
%
%as provided by formulae 6.1.29 and 6.1.30 of \citet{abramowitz}.

\smallskip

The same procedure applies for the calculation of the denominator in \eqref{eq:expbnd}. Other expressions of $b_{n,d}$ have been provided by \citet{artggipcu18} and \citet{lantuejoul2019}, but they are valid only when $d=2$. 

%The explicit computation of the proposed coefficients can be performed using the induction formula 
%
%$$ b_{n+2,d} = \frac{\lambda + n + 2}{\lambda + n} \, \frac{n^2 + \nu^2}{4 (\lambda + 1 +n)^2 + \nu^2} \, b_{n,d}, $$
%
%starting from 
%
%\begin{align*}
%b_{0,d} &= \exp ( - \pi \nu / 2) \, \frac{ \lvert \Gamma (\lambda + 1) \rvert^2}{ \displaystyle \left\lvert \Gamma 
%\left( \lambda +1 +  i \nu / 2 \right) \right\rvert^2} \\
%b_{1,d} &=  \frac{\nu}{2} \,  \exp ( - \pi \nu / 2) \, \frac{ \Gamma (\lambda) \, \Gamma (\lambda + 2)}{ \displaystyle 
%\left\lvert \Gamma \left( \lambda + 3/2 +  i \nu / 2 \right) \right\rvert^2}.  
%\end{align*}
%
%Note that, if $d$ is odd, the computation of the Schoenberg coefficient \eqref{eq:expbnd} is straightforward: 
%
%$$ b_{n,2k+1} = \frac{C(\nu,n) \, (k+n) \, (k-1)! \, k! \, 4^{k+1}}{\displaystyle \prod_{\ell=0}^k \bigl( ( 2 \ell +n )^2 +  %\nu^2 \bigr)}. $$
%

%

%%%%%%%%%%%%%%%%%%%%%%%%%%%%%%%%%%%%%%%%%%%%%%%%%%%%%%%%%%%%%%%%%%%%%%%%%%%%%%%%%%%%%%%%%%%%%%%%%

%%%%%%%%%%%%%%%%%%%%%%%%%%%%%%%%%%%%%%%%%%%%%%%%%%%%%%%%%%%%%%%%%%%%%%%%%%%%%%%%%%%%%%%%%%%%%%%%%

%\section*{\refname}
\bibliographystyle{apalike}
%\biboptions{authoryear}
\bibliography{mybib}

\end{document}